\documentclass[12pt]{article}
\usepackage{amsfonts,amssymb,latexsym}

\newtheorem{theorem}{Theorem}[section]
\newtheorem{lemma}[theorem]{Lemma}

\newtheorem{proposition}[theorem]{Proposition}

\newenvironment{proof}
{\par\addvspace{0.3cm}\noindent{\rm Proof. }}
{\nopagebreak\mbox{}\hfill $\Box$\par\addvspace{0.25cm}}
\newenvironment{proofof}[1]
{\par\addvspace{0.3cm}\noindent{\rm Proof of #1. }}
{\nopagebreak\mbox{}\hfill $\Box$\par\addvspace{0.25cm}}
\renewcommand{\Re}{\mathrm{Re\,}}
\renewcommand{\Im}{\mathrm{Im\,}}
\newcommand{\sech}{\mathrm{sech\,}}
\newcommand{\R}{{\mathbb R}}
\newcommand{\C}{{\mathbb C}}
\newcommand{\Z}{{\mathbb Z}}
\newcommand{\T}{{\mathbb T}}
\newcommand{\W}{\mathcal{W}}
\newcommand{\cF}{\mathcal{F}}
\newcommand{\cP}{\mathcal{P}}
\newcommand{\cL}{\mathcal{L}}

\newcommand{\Li}{L^\iy(\T)}

\renewcommand{\kappa}{\varkappa}
\renewcommand{\rho}{\varrho}

\newcommand{\be}{\begin{equation}}
\newcommand{\ee}{\end{equation}}
\newcommand{\ds}{\displaystyle}
\newcommand{\ts}{\textstyle}
\newcommand{\bqn}{\begin{eqnarray}}
\newcommand{\eqn}{\end{eqnarray}}
\newcommand{\nn}{\nonumber}
\newcommand{\ba}{\begin{array}}
\newcommand{\ea}{\end{array}}

\newcommand{\al}{\alpha}
\newcommand{\iv}{^{-1}}
\newcommand{\iy}{\infty}
\newcommand{\eps}{{\varepsilon}}

\newcommand{\ovl}{\overline}
\newcommand{\trace}{\mathrm{trace}}

\newcommand{\twomat}[1]{\left(\ba{cc}#1\ea\right)}

\newcommand{\ta}{\tilde{a}}
\newcommand{\tb}{\tilde{b}}

\newcommand{\el}{\ell^2}

\newcommand{\hv}{\hat{v}}
\newcommand{\hu}{\hat{u}}
\newcommand{\ha}{\hat{a}}
\newcommand{\hK}{\hat{K}}

\newcommand{\Pe}{\Pi_{[\eps,1]}}
\newcommand{\Per}{\Pi_{[\eps,\sqrt{\eps}]}}
\renewcommand{\Pr}{\Pi_{[\sqrt{\eps},1]}}
\newcommand{\Pei}{\Pi_{[1,\eps\iv]}}

\begin{document}

\date{}
\title{On the asymptotics of certain Wiener-Hopf-plus-Hankel determinants}
\author{Estelle L. Basor
	\thanks{ebasor@calpoly.edu, Supported in part by NSF Grant DMS-0200167.}\\
	Department of Mathematics\\
	California Polytechnic State University\\
	San Luis Obispo, CA 93407, USA\\ \\
         Torsten Ehrhardt
	\thanks{ehrhardt@math.ucsc.edu}\\
	Department of Mathematics\\
         University of California\\
         Santa Cruz, CA 95064, USA} 	
\maketitle

\begin{abstract}
In this paper we determine the asymptotics of the determinants
of truncated Wiener-Hopf plus Hankel operators
$\det (W_R(a)\pm H_R(a))$ as $R\to\iy$ for symbols
$a(x)=( x^2/(1+x^2))^\beta$ with the parameter $\beta$ being of small
size.
\end{abstract}




\section{Introduction}

For a function $a$ defined on the real line $\R$ such that 
$a-1\in L^1(\R)$ the truncated Wiener-Hopf and Hankel operators 
acting on $L^2[0,R]$ with symbol $a$ are defined by
\bqn
W_R(a) &:& f(x)\mapsto g(x)=f(x)+\int_{0}^R k(x-y)f(y)\,dy,\\[1ex]
H_R(a) &:& f(x)\mapsto g(x)=\int_{0}^R k(x+y)f(y)\,dy,
\eqn
where $k$ is the Fourier transform of $a-1$,
\bqn\label{f.FT}
k(x)=\frac{1}{2\pi}\int_{-\iy}^\iy (a(\xi)-1)e^{-i\xi x}\,d \xi.
\eqn
It is well known that under the above assumption
the operators $W_R(a)-I$ and $H_R(a)$ are trace class operators.
Hence the determinants
$$
\det(W_R(a)\pm H_R(a))
$$
are well-defined.

The purpose of this paper is to determine the asymptotics of these determinants
as $R\to\iy$ for a particular class of even generating functions which have a
single singularity at $x=0$. Before explaining the scope of this paper in
more detail, let us briefly review related problems.

The asymptotics of Wiener-Hopf determinants $\det W_R(a)$ as $R\to\iy$
for sufficiently smooth nonvanishing functions $a$ with winding number zero
are described by the Akhiezer-Kac formula (see, e.g. \cite{BS} and the references
therein). A more complicated situation occurs when the symbol
$a$ possesses singularities such as jumps, zeros, or poles.
Let $\hu_\beta$ and $\hv_\beta$ be the functions
\be\label{f.vuh}
\hv_{\beta}(x):=\left(\frac{x^2}{x^2+1}\right)^{\beta},\qquad
\hu_{\beta}(x):=\left(\frac{x-0 i}{x-i}\right)^{-\beta}
\left(\frac{x+0 i}{x+i}\right)^{\beta}.
\ee
Notice that $\hv_\beta$ has a zero or a pole at $x=0$, while
$\hu_\beta$ has a jump discontinuity at $x=0$ whose size is determined
by the parameter $\beta$. If the symbol is of the Fisher-Hartwig form,
\bqn
\ha(x)&=& b(x)\prod_{r=1}^R \hv_{\alpha_r}(x-x_r)\hu_{\beta_r}(x-x_r),
\eqn
where $|\Re\alpha_r|<1/2$, $|\Re\beta_r|<1/2$, $x_1,\dots,x_R\in\R$ are
distinct, and $b$ is a sufficiently smooth function satisfying the assumptions
of the Akhiezer-Kac formula, then the asymptotics of the
determinants are described by the continuous analogue of the Fisher-Hartwig 
conjecture. One minor complication is encountered since, except in special cases, 
the symbol $\ha$ does not belongs to $L^1(\R)$, but only to $L^2(\R)$. 
Because then the above operators are only Hilbert-Schmidt one has to consider 
regularized determinants $\det_2(I+K)=\det(I+K)e^{-K}$. 
The asymptotic formula for such Wiener-Hopf determinants reads
\bqn\label{f.aW}
{\ts \det_2} W_R(\ha)&\sim& G_2[\ha]^R\, R^\Omega E,\qquad R\to\iy,
\eqn
where $\Omega=\sum_{r=1}^R(\alpha_R^2-\beta_R^2)$, $G_2[\ha]$ is a regularized version
of the geometric means of $\ha$, and $E$ is a complicated constant. 

Formula (\ref{f.aW}) has not yet been proved general (see \cite{BSW} for the proof in a special case
where $\al_r=0$ for all $r$), 
but it is very likely that such a proof can be accomplished with the help of two main ingredients.
One of these is a localization theorem for Wiener-Hopf determinants, which had to be
analogous to a corresponding (well-known) localization theorem for Toeplitz 
determinants with Fisher-Hartwig symbols \cite{BS}. (The outline of a possible proof of such a theorem
has been communicated to us by A. B\"ottcher, but the details still need to be verified.)

The localization reduces the problem to symbols that are ``pure'' Fisher-Hartwig symbols, that is where $R=1$ and $b(x) \equiv 1.$ 
This last problem was outstanding for a long time and was recently solved
by one of the authors and Widom \cite{BW1}. They made use of the so-called
Borodin-Okounkov formula \cite{BO} (see also \cite{BW0,Bot}) to compare the
asymptotics of $\det_2 W_R(\ha)$ with the (known) asymptotics
of a Toeplitz determinant $\det T_n(a)$ where $R\sim 2n$ and $n,R\to\iy$.
The Borodin-Okounkov identity is an exact identity for both the Toeplitz and Wiener-Hopf determinants and made the comparisons possible.

We will do something very similar in this paper, in the sense that we will also make a comparison to already known asymptotics. These will involve the discrete analogue of the sum of the finite Wiener-Hopf and Hankel operators.
The discrete analogues  are the
Toeplitz and Hankel matrices,
\be
T_n(a)=(a_{j-k})_{j,k=0}^{n-1},\qquad
H_n(a)=(a_{j+k+1})_{j,k=0}^{n-1}.
\ee
Here $a\in L^1(\T)$ is a function defined on the unit circle
$\T=\{t\in\C\;:\;|t|=1\}$ with Fourier coefficients
$$
a_k =\frac{1}{2\pi}\int_0^{2\pi}a(e^{i\theta})e^{-ik\theta}\, d\theta,\qquad
k\in\Z.
$$
The asymptotics of Toeplitz determinants have a long and interesting history.
For the latest results and more information we refer to \cite{E}.

The study of the asymptotics of Toeplitz-plus-Hankel determinants
$\det(T_n(a)\pm H_n(a))$ was begun recently. 
The main interest is in even symbols (i.e., $a(t)=a(t\iv)$, $t\in\T$).
In this case the results for Fisher-Hartwig type symbols
are nearly complete \cite{BE1}. Some results have been obtained also 
for non-even symbols \cite{BE0}.

Let us now return to the topic of this paper, namely the asymptotics of
Wiener-Hopf-plus-Hankel determinants, 
$$
\det(W_R(\ha)\pm H_R(\ha)).
$$
First of all, the case of smooth, nonvanishing and even functions 
(i.e, $\ha(x)=\ha(-x)$, $x\in\R$) follows from (more general) results in \cite{BE2}. 
In regard to Fisher-Hartwig type symbols only the case
of a function $\ha(x)=\hu_{\beta}(x-1)\hu_{-\beta}(x+1)$ 
(which is a even piecewise constant function with two jump discontinuities)
was treated recently in \cite{BEW}. 

In this paper we consider the case of a function $\ha(x)=\hv_\beta(x)$, 
which is an even function having a zero or a pole at $x=0$. 
In order to state the main result  we introduce
\be
\hat{D}_R^+(\beta):=\det\Big[W_R(\hv_\beta)+ H_R(\hv_\beta)\Big],
\qquad
\hat{D}_R^-(\beta):=\det\Big[W_R(\hv_\beta)- H_R(\hv_\beta)\Big].
\ee

The natural assumption on $\beta$ is that $\Re\beta>-1/2$ since then 
$\hv_\beta\in L^1(\R)$. Moreover, because $W_R(\hv_\beta)$ and $H_R(\hv_\beta)$ 
are analytic operator valued functions with respect to $\beta$, 
the functions $\hat{D}^\pm_R(\beta)$ are analytic on the set 
of all $\beta\in\C$ for which $\Re\beta>-1/2$.

\begin{theorem}\label{t.main}\
\begin{itemize}
\item[(a)]
If $-1/2<\Re\beta<3/2$, then
\bqn
\hat{D}_R^+(\beta)
&\sim&
e^{-\beta R}
R^{\beta^2/2-\beta/2} (2\pi)^{\beta/2} 2^{-\beta^2+\beta/2}
\frac{G(1/2)}{G(1/2+\beta)},\quad R\to\iy.\label{f.DR+}
\eqn
\item[(b)]
The function $\hat{D}^-_R(\beta)$ admits an analytic continuation
onto the set of all $\beta\in\C$ for which $\Re\beta>-3/2$.
Moreover, if $-1<\Re\beta<1/2$, then 
\bqn
\hat{D}_R^-(\beta)
&\sim&
e^{-\beta R}
R^{\beta^2/2+\beta/2} (2\pi)^{\beta/2} 2^{-\beta^2-\beta/2}
\frac{G(3/2)}{G(3/2+\beta)},\quad R\to\iy.\label{f.DR-}
\eqn
\end{itemize}
\end{theorem}

Therein $G(z)$ is the Barnes $G$-function \cite{Bar},
which is an entire function defined by
\bqn
G(1+z) &=& (2\pi)^{z/2} e^{-(z+1)z/2-\gamma_Ez^2/2}
\prod_{k=1}^\iy\Big((1+z/k)^ke^{-z+z^2/(2k)}\Big)
\eqn
with $\gamma_E$ equal to Euler's constant. Notice that the Barnes function
has the remarkable property that $G(1+z)=\Gamma(z) G(z)$.

The proof of the above theorem will be given in Section \ref{s3.6}.
However, the first statement in part (b) concerning the analytic continuablity  
will already follow from Proposition \ref{p2.13} in Section \ref{sec:3:4}.

The assumption $-1<\Re\beta<1/2$ rather than $-3/2<\Re\beta<1/2$
in part (b) seems to be too restrictive. Unfortunately, 
we have not been able to remove it.

It is interesting to observe that Theorem \ref{t.main} implies the asymptotics
for Wiener-Hopf determinants
\bqn
\det W_{2R}(\hv_\beta) &\sim& e^{-2\beta R} R^{\beta^2}
\frac{G(1+\beta)^2}{G(1+2\beta)},
\qquad R\to\iy,
\eqn
which has been proved in \cite{BW1}. To see this one has to use
the formula
\bqn\label{W2R}
\det W_{2R}(\hv_\beta)&=&\det\Big[W_R(\hv_\beta)+H_R(\hv_\beta)\Big]
\cdot\det\Big[W_R(\hv_\beta)-H_R(\hv_\beta)\Big],
\eqn
which can be easily proved by observing that $W_{2R}(\hv_{\beta})$ 
can be identified with the block operator
$$
\left(\ba{cc} W_{R}(\hv_{\beta})&H_{R}(\hv_{\beta})\\
H_{R}(\hv_{\beta})&W_{R}(\hv_{\beta})\ea\right).
$$
One has also to make use of a consequence of the 
duplication formula for the Barnes $G$-function, which will be stated 
below in (\ref{f.dupl1}).

An outline of the paper is as  follows. In the next section, we review the asymptotics in the discrete case. The final section (Section 3) is divided into several subsections and contains the proof of Theorem
\ref{t.main}. We first (Sect.~3.1 and 3.2) review the basic operator theory facts that are needed and proceed then (Sect.~3.3 and 3.4) with identifying the determinants $\hat{D}_n^\pm (\beta)$ and their
discrete analogues with different kinds of determinants. 
In Section 3.5 we prove some theorems concerning the strong  and trace class convergence of certain operators that naturally occur in our proof. Lemma \ref{l2.13} is really the basic formula that allows us to compare the desired determinants with those of the discrete analogues. Finally, in Section 3.6, we
will complete the proof of the main results, and,  for completeness sake, we will  compute some Fredholm determinants that occur that may be interesting in other settings. 

We end this introduction by pointing out that the asymptotics of the determinants of Wiener-Hopf-plus-Hankel operators with the type of discontinuity considered have important applications. The computation of such asymptotics is a crucial step in the work of the second author's recent proof of the complete asymptotics of the Fredholm determinant of the sine kernel \cite{E2} (see also \cite{D1,W3}). It is also true that such operators occur in the Laguerre random matrix ensemble in a very natural way when one considers special parameters ($\nu = \pm 1/2$ for Bessel operators; see \cite{BE2}). Computing the asymptotics for singular symbols yields information about certain discontinuous random variables.


\section{The asymptotics in the discrete case} 

In what follows we are going to recall the results
about the asymptotics of the determinants
\be
D_n^+(\beta):=\det\Big[T_n(v_\beta)+H_n(v_\beta)\Big],\qquad
D_n^-(\beta):=\det\Big[T_n(v_\beta)-H_n(v_\beta)\Big],
\ee
as $n\to\iy$, which were established in \cite{BE0}.
Therein $v_\beta$ is the function 
\bqn
\label{f1.v}
v_\beta(e^{i\theta}) &:=& (2-2\cos\theta)^{\beta},
\qquad\Re\beta>-1/2.
\eqn
Let us also introduce the function 
\bqn
\label{f1.u}
u_{\beta}(e^{i\theta}) &:=& e^{i\beta(\theta-\pi)},
\qquad\qquad
0<\theta<2\pi.
\eqn
It is easily seen that $D_n^\pm(\beta)$ are analytic in $\beta$
for $\Re\beta>-1/2$.  In the following theorem we will provide some information about the
analytic continuability of $D_n^\pm(\beta)$ with respect to $\beta$ and
about the asymptotics of $D_n^\pm(\beta)$ as $n\to\iy$ for fixed $\beta$.

\begin{theorem}\label{t.Ddiscrete}\
\begin{itemize}
\item[(a)] For each $n\ge1$ the function $D_n^+(\beta)$ is analytic
in $\beta$ on $U_+:=\C\setminus\{-1/2,-3/2,-5/2,\dots\}$.
Moreover, for $\beta\in U_+$,
\bqn
D_n^+(\beta) &\sim&
n^{\beta^2/2-\beta/2}\,(2\pi)^{\beta/2}\,2^{-\beta^2/2}\,
\frac{G(1/2)}{G(1/2+\beta)},\qquad n\to\iy.
\eqn
\item[(b)] For each $n\ge1$ the function $D_n^-(\beta)$ is analytic
in $\beta$ on $U_-=\C\setminus\{-3/2,-5/2,-7/2,\dots\}$.
Moreover, for $\beta\in U_-$, 
\bqn
D_n^-(\beta) &\sim&
n^{\beta^2/2+\beta/2}(2\pi)^{\beta/2}\,2^{-\beta^2/2}\,
\frac{G(3/2)}{G(3/2+\beta)},\qquad n\to\iy.
\eqn
\end{itemize}
\end{theorem}
\begin{proof}
{}From the proof of Thm.~7.7 in \cite{BE0}, it follows that
\bqn
\frac{\det\Big[T_n(v_\beta)\pm H_n(v_\beta)\Big]}
{\det\Big[T_n(u_{-\beta})\mp H_n(u_{-\beta})\Big]}
&=& \prod_{k=0}^{n-1}
\frac{\Gamma(1+2\beta+k)\Gamma(1-\beta+k)}
{k!\,\Gamma(1+\beta+k)}\nn\\[1ex]
&=&
\frac{G(1+2\beta+n)G(1-\beta+n)}{G(1+n)G(1+\beta+n)}
\frac{G(1+\beta)}{G(1+2\beta)G(1-\beta)}.\nn
\eqn
(Notice the different meaning of the notation of $u_\beta$ used there.)
Furthermore, the proof of Thms.~6.2 and 6.3 in \cite{BE0} implies that
$$
\ba{l}
\det\Big[T_n(u_\alpha)+H_n(u_\alpha)\Big] 
\quad=
\nn\\[3ex]\ds
\quad
(2\pi)^{\al/2}\, 2^{\al^2/2+1}\,
\frac{G(1/2-\al)G(1+\al)G(1-\al)}{G(1/2)}
\frac{G(2n)G(2n-2\al)}{G(2n+1-\al)G(2n-1-\al)}
\nn\\[3ex]\ds
\quad
\times
\frac{G(n+3/2-\al)G(n+1)G(n+1-\al)G(n-1/2-\al/2)G(n-\al/2)^2G(n+1/2-\al/2)}
{G(n+1/2-\al)^2G(n+1/2)G(n)G(n+1-2\al)G(n-\al)G(n+\al+1)}\nn
\ea
$$
and
$$
\ba{l}
\det\Big[T_n(u_\alpha)-H_n(u_\alpha)\Big] 
\quad=
\nn\\[3ex]\ds
\quad
(2\pi)^{\al/2}\, 2^{\al^2/2+1}\,
\frac{G(3/2-\al)G(1+\al)G(1-\al)}{G(3/2)}
\frac{G(2n)G(2n-2\al)}{G(2n+1-\al)G(2n-1-\al)}
\nn\\[3ex]\ds
\quad
\times
\frac{G(n+3/2)G(n+1)G(n+1-\al)G(n-1/2-\al/2)G(n-\al/2)^2G(n+1/2-\al/2)}
{G(n+1/2-\al)G(n+1/2)^2G(n)G(n+1-2\al)G(n-\al)G(1+\al+n)}.\nn
\ea
$$
Combining these results we obtain
\bqn
D_n^+(\beta)
&=&
(2\pi)^{-\beta/2}\, 2^{\beta^2/2+1}\,
\frac{G(3/2+\beta)G(1+\beta)^2}{G(3/2)G(1+2\beta)}
\frac{G(2n)G(2n+2\beta)}{G(2n+1+\beta)G(2n-1+\beta)}
\nn\\[1ex]
&&
\times\,
\frac{G(n+3/2)G(n-1/2+\beta/2)G(n+\beta/2)^2G(n+1/2+\beta/2)}
{G(n)G(n+1/2)^2G(n+1/2+\beta)G(n+\beta)}\nn
\eqn
and
\bqn
D_n^-(\beta)
&=&
(2\pi)^{-\beta/2}\, 2^{\beta^2/2+1}\,
\frac{G(1/2+\beta)G(1+\beta)^2}{G(1/2)G(1+2\beta)}
\frac{G(2n)G(2n+2\beta)}{G(2n+1+\beta)G(2n-1+\beta)}
\nn\\[1ex]
&&
\times\,
\frac{G(n+3/2+\beta)G(n-1/2+\beta/2)G(n+\beta/2)^2G(n+1/2+\beta/2)}
{G(n) G(n+1/2)G(n+1/2+\beta)^2G(n+\beta)}.\nn
\eqn
Using the duplication formula for the $G$-function \cite[p.~291]{Bar},
\bqn\label{f.dupl0}
G(z)G(z+1/2)^2G(z+1) &=& G(1/2)^2\pi^z2^{-2z^2+3z-1}G(2z),
\eqn
it now follows that
\bqn\label{f.dupl1}
\frac{G(1/2+\beta)G(1+\beta)^2G(3/2+\beta)}{G(1+2\beta)}
= (2\pi)^{\beta}2^{-2\beta^2}G(1/2)G(3/2)
\eqn
and 
\bqn
&&
\frac{G(2n)G(2n+2\beta)}{G(2n+1+\beta)G(2n-1+\beta)}
\quad=\quad 2^{\beta^2-1}\,\times\nn\\
&& \quad
\frac{G(n)G(n+1/2)^2G(n+1)G(n+\beta)G(n+1/2+\beta)^2G(n+1+\beta)}
{G(n-1/2+\beta/2)G(n+\beta/2)^2G(n+1/2+\beta/2)^2G(n+1+\beta/2)^2G(n+3/2+\beta/2)}
.\nn
\eqn
Hence
\bqn
D_n^+(\beta) &=&
(2\pi)^{\beta/2}\,2^{-\beta^2/2}\,
\frac{G(1/2)}{G(1/2+\beta)}
\nn\\[1ex]
&&\times
\frac{G(n+3/2)G(n+1)G(n+1+\beta)G(n+1/2+\beta)}
{G(n+1/2+\beta/2)G(n+1+\beta/2)^2G(n+3/2+\beta/2)}\nn
\eqn
and
\bqn
D_n^-(\beta) &=&
(2\pi)^{\beta/2}\,2^{-\beta^2/2}\,
\frac{G(3/2)}{G(3/2+\beta)}
\nn\\[3ex]
&&\times
\frac{G(n+1/2)G(n+1)G(n+1+\beta)G(n+3/2+\beta)}
{G(n+1/2+\beta/2)G(n+1+\beta/2)^2G(n+3/2+\beta/2)}.\nn
\eqn
{}From these identities it can be concluded that $D_n^\pm(\beta)$ can
continued analytically to all of $U_\pm$. 
Notice that the zeros in the denominator cancel with the zeros of the term 
$G(n+1+\beta)$ in the numerator.

In order to obtain the asymtotic result, we apply the formula
\bqn\label{f.BarAsy}
\prod_{r=1}^R\frac{G(1+x_r+n)}{G(1+y_r+n)}&\sim& n^{\omega/2},\qquad n\to\iy,
\eqn
which holds under the assumption $x_1+\dots+x_R=y_1+\dots+y_R$ and 
with  the constant $\omega=x_1^2+\dots+x_R^2-y_1^2-\dots-y_R^2$.
This asymptotic formula has been proved, e.g., in Lemma 6.1 of \cite{BE0}.
\end{proof}

Once again it is interesting to remark that from the asymptotics established
in the previous theorem and from the identity
\bqn\label{det.Tn}
\det T_{2n}(v_\beta) &=& D_n^+(\beta)D_n^-(\beta)
\eqn
the well-known asymptotics
\bqn
\det T_{2n}(v_\beta)&\sim& (2n)^{\beta^2}\frac{G(1+\beta)^2}{G(1+2\beta)}
\eqn
follow by using the consequence (\ref{f.dupl1}) of the duplication formula for the 
Barnes function.  As in the case of (\ref{W2R}) formula (\ref{det.Tn}) can be 
proved by identifying the symmetric matrix $T_{2n}(v_\beta)$ with a 
two-by-two block matrix having  the entries $T_{n}(v_{\beta})$ and
$H_{n}(v_{\beta})$.


\section{Proof of the main results}

\subsection{Preliminary facts}

An operator $A$ acting on a Hilbert space $H$ is called a trace class operator 
if it is compact and if the series consisting of the singular values $s_n(A)$
(i.e., the eigenvalues of $(A^*A)^{1/2}$ taking multiplicities into account)
converges. The norm
\bqn
\|A\|_{1} &=& \sum_{n\ge1} s_n(A)
\eqn
makes the set of all trace class operators into a Banach space, which forms
also a two-sided ideal in the algebra of all linear bounded operators on $H$.
Moreover, the estimates
$\|AB\|_1\le\|A\|_1\|B\|$ and $\|BA\|_1\le\|A\|_1\|B\|$ hold,
where $A$ is a trace class operator and $B$ is a bounded operator with the
operator norm $\|B\|$.

A useful property is that  if $B$ is a trace class operator, if
$A_n\to A$ converges strongly on $H$ and if $C_n^*\to C^*$ converges
strongly on $H$, then $A_nB C_n\to ABC$ in the trace norm. Therein
$C^*$ stands for the Hilbert space adjoint of the operator $C$.

If $A$ is a trace class operator, then the operator trace
``$\mathrm{trace}(A)$'' and the operator determinant ``$\det (I+A)$'' are
well defined. For more information concerning these concepts we refer to
\cite{GK}.

A sequence of operators $A_n$ defined on a Hilbert space $H$ is called stable
if the operators $A_n$ are invertible for all sufficiently large $n$ and
if 
$$
\sup\limits_{n\ge n_0}\|A_n\iv\|_{\cL(H)}\;<\;\iy.
$$
A sequence of linear bounded operators $A_n$
on a Hilbert space $H$ is said to converge strongly on $H$ to an operator $A$ if
$A_n x\to Ax$ for all $x\in H$.

\begin{lemma}\label{l3.2}
Let $A_n$ be a sequence of linear bounded operators on a Hilbert space $H$
such that $A_n\to A$ strongly. Then $A_n\iv\to A\iv$ strongly if and only if
the sequence $A_n$ is stable.
\end{lemma}
\begin{proof}
The ``if'' part can be proved by using the estimate
$$
\|A_n\iv x-A\iv x\|\le \|A_n\iv\|\cdot \|(A-A_n)A\iv x\|.
$$
The ``only if'' part of the lemma follows from the Banach-Steinhaus Theorem.
\end{proof}

In what follows we consider some concrete classes of
linear bounded operators.

For a function $a\in \Li$ with Fourier coefficients $\{a_n\}_{n=-\iy}^\iy$,
the Toeplitz and Hankel operators are linear bounded operators
acting on $\el=\el(\Z_+)$ defined by the infinite matrices
\be\label{f.THdef}
T(a)=(a_{j-k})_{j,k=0}^\iy,\qquad
H(a)=(a_{j+k+1})_{j,k=0}^\iy.
\ee
The connection to $n\times n$ Toeplitz and Hankel matrices is given by
\bqn
P_n T(a)P_n \cong T_n(a),\qquad
P_n H(a) P_n \cong H_n(a),
\eqn
where $P_n$ is the finite rank projection operator on $\el$
\bqn
P_n: (x_0,x_1,\dots)\in \el\mapsto (x_0,\dots,x_{n-1},0,\dots)\in\el.
\eqn
Toeplitz and Hankel operators satisfy the following well-known formulas,
\bqn
T(ab)&=& T(a)T(b)+H(a)H(\tb),\label{f.Tab}\\[1ex]
H(ab)&=& T(a)H(b)+H(a)T(\tb),\label{f.Hab}
\eqn
where $\tb(t):=b(t\iv)$, $t\in\T$.

For a functions $a\in L^\iy(\R)$ the Wiener-Hopf operator and 
the Hankel operator acting on $L^2(\R_+)$ are defined by
\bqn\label{f.WH}
W(a)&=& P_+  \cF M(a)\cF\iv P_+|_{L^2(\R_+)},
\\[1ex]
H(a)&=&P_+  \cF M(a)\cF\iv \hat{J}P_+|_{L^2(\R_+)},
\eqn
where $\cF$ stands for the Fourier transform acting on $L^2(\R)$,
$M(a)$ stands for the multiplication operator on $L^2(\R)$, $P_+=M(\chi_{\R_+})$,
and $(\hat{J}f)(x)=f(-x)$.
If $a\in L^1(\R)\cap L^\iy(\R)$, then $W(a)$ and $H(a)$ are integal operators 
on $L^2(\R)$ with the kernel $k(x-y)$ and $k(x+y)$, respectively, where
$k(x)$ is the Fourier transform (\ref{f.FT}) of $a$. We remark that 
\bqn
W(ab) &=& W(a)W(b)+H(a)H(\tb),
\\[1ex]
H(ab) &=& W(a)H(b)+H(a)W(\tb),
\eqn
where $\tilde{b}(x):=b(-x)$, $x\in\R$. Moreover,
\be
W_R(a)= P_R W(a) P_R|_{L^2[0,R]},\qquad H_R(a)= P_R H(a) P_R|_{L^2[0,R]},
\ee
where $P_R=M(\chi_{[0,R]})$.

It is important to note that Wiener-Hopf and Hankel operators are
related to their discrete analogues by a unitary transform
$S:\el\to L^2(\R_+)$,
\bqn\label{f.cd}
T(a)=S^* W(\ha) S,\qquad
H(a)= S^* H(\ha)S,
\eqn
where the symbols are related by
\bqn\label{f.frac}
\ha(x)=a\left(\frac{1+ix}{1-ix}\right).
\eqn
(The use of the same notation for the continuous and the discrete Hankel
operators should not cause confusion.) 
For sake of further reference, let us introduce the mapping
\be\label{f.Phi}
\Phi:A\in \cL(L^2(\R_+))\mapsto S^* A S \in \cL(\el).
\ee
The unitary transform $S$ is given explicitly by the composition
$S=\cF U \cF_{\mathrm{d}}\iv$, where
$$
\el\stackrel{\cF_{\mathrm{d}}\iv}{\longrightarrow}
H^2(\T)\stackrel{U}{\longrightarrow}
H^2(\R)\stackrel{\cF}{\longrightarrow} L^2(\R_+).
$$
Therein $H^2(\T)$ and $H^2(\R)$ are the Hardy spaces with respect to $\T$ and
$\R$,
\bqn
H^2(\T)&=&\Big\{\;f\in L^2(\T)\;:\;f_k=0\mbox{ for all } k<0\;\Big\},
\nn\\[1ex]
H^2(\R)&=&\Big\{\;f\in L^2(\R)\;:\;(\cF f)(x)=0\mbox{ for all } x<0\;\Big\},
\eqn
$\cF_{\mathrm{d}}\iv:\{x_{n}\}_{n=0}^\iy\mapsto f(t)=\sum_{n=0}^\iy x_n t^n$ is
the inverse discrete Fourier transform, and 
$$
(Uf)(x)=\frac{1}{\sqrt{\pi}(1-ix)}f\left(\frac{1+ix}{1-ix}\right).
$$

Under the unitary transform $\cF_{\mathrm{d}}: H^2(\T)\to \el $, the Toeplitz and Hankel operators
can be identified with operators acting on $H^2(\T)$,
\be\label{f.THid}
T(a)\cong PM(a)P|_{H^2(\T)},\qquad
H(a)\cong PM(a)JP|_{H^2(\T)},
\ee
where  $P$ is the Riesz projection on $L^2(\T)$, $M(a)$ is
the multiplication operator on $L^2(\T)$, and $(Jf)(t)=t\iv f(t)$,
$t\in\T$.

A sequence of functions $a_n\in \Li$ is said to converge to $a\in \Li$ in measure
if for each $\eps>0$ the Lebesgue measure of the set
$$
\Big\{\;t\in\T\;:\; |a_n(t)-a(t)|\ge \eps\;\Big\}
$$
converges to zero.

\begin{lemma}\label{l3.1}
Assume that $a_n\in \Li$ are uniformly bounded and converge to $a\in \Li$
in measure. Then
$$
T(a_n)\to T(a)\quad\mbox{ and } \quad
H(a_n) \to H(a)
$$
strongly on $\el$, and the same holds for the adjoints.
\end{lemma}
\begin{proof}
We use the identification (\ref{f.THid}).
If $a_n$ converges in measure to $a$ and is uniformly bounded, then
$a_n$ also converges to $a$ in the $L^2$-norm. Hence for all
$f\in L^\iy$, we have $a_n f\to af$ in the $L^2$-norm. Using an approximation
argument and the uniform boundedness of $a_n$, it follows
that $M(a_n)\to M(a)$ strongly on $L^2(\T)$.
Hence the corresponding Toeplitz and Hankel operators
converge strongly on $H^2(\T)$, too. Since
$T(a_n)^*=T(a_n^*)$ and $H(a_n)^*=H(\tilde{a}_n^*)$, this
holds also for the adjoints.
\end{proof}


\subsection{Invertibility of operators $I\pm H(u_{\beta})$}

In this section we prove that operators  of the form $I\pm H(u_{\beta})$ 
are invertible for certain $\beta$. We think of the Hankel operators 
as discrete ones acting on $\el$ (or, equivalently, on $H^2(\T)$). Obviously, these 
invertibility results can be extended with the help of (\ref{f.frac}) and (\ref{f.Phi})
to operators $I\pm H(\hu_{\beta})$ where continuous Hankel 
operators acting on $L^{2}(\R_{+})$ are involved.

For $\tau\in\T$ and $\beta\in\C$ we introduce the functions
\bqn
\eta_{\beta}(t)=(1-t)^{\beta},\qquad
\xi_{\beta}(t)=(1-1/t)^{\beta},
\eqn
where these functions are analytic in an open neighborhood of
$\{\;z\in\C\;:\;|z|\le 1,\;z\neq 1\}$ and
$\{\;z\in\C\;:\;|z|\ge 1,\;z\neq 1\}\cup\{\iy\}$, resp.,
and the branch of the power function is chosen in such a way that
$\eta_{\beta}(0)=1$ and $\xi_{\beta}(\iy)=1$.
Notice that
\bqn\label{f.uv-rel}
v_{\beta}(t)=\eta_{\beta}(t)\xi_{\beta}(t),\qquad
u_{\beta}(t)=\eta_{\beta}(t)\xi_{-\beta}(t),\qquad
u_{\beta+n}(t)=(-t)^n u_{\beta}(t).
\eqn

The essential spectrum $\mathrm{sp}_{\mathrm{ess}} A$
of a linear bounded operator $A$ defined on a
Hilbert space is the set of all $\lambda \in \C$ for which
$A-\lambda I$ is not a Fredholm operator.

\begin{proposition}\label{p2.1}
Let $\beta\in\C$. Then
\begin{itemize}
\item[(a)]
$I+H(u_{\beta})$ is Fredholm on $\el$ if and only if
$\Re \beta\notin \frac{1}{2}+2\Z$,
\item[(b)]
$I-H(u_{\beta})$ is Fredholm on $\el$ if and only if
$\Re \beta\notin -\frac{1}{2}+2\Z$.
\end{itemize}
\end{proposition}
\begin{proof}
We use a result of Power \cite{Po1} in order to determine the
essential spectrum of the Hankel operator operators $H(u_\beta)$.
It says that the essential spectrum is a union of intervals in the complex
plane, namely
\bqn\label{f.Spess}
\mathrm{sp}_{\mathrm{ess}}H(b) &=&[0,ib_{-1}]\cup[0,-ib_1]\cup
\bigcup_{\textstyle{\tau\in\T\atop \mathrm{Im} \tau>0}}
\left[-i\sqrt{b_\tau b_{\bar{\tau}}},i \sqrt{b_\tau b_{\bar{\tau}}} \right].
\eqn
Therein we use the notation $b_\tau=(b(\tau+0)-b(\tau-0))/2$ with
$b(\tau\pm0)=\lim_{\eps\to+0} b(\tau e^{i\eps})$.
This result can also be obtained from the
more general results contained in \cite{Po2} and \cite[Sect.~4.95--102]{BS}.

Clearly, $b_\tau=0$ for $b=u_\beta$ if $\tau\neq1$. In the case $\tau=1$
we have $b_1=-i\sin(\beta\pi)$. Hence
$$
\mathrm{sp}_{\mathrm{ess}}H(u_{\beta})=[0,-\sin(\pi\beta)],
$$
from which the assertion is easy to conclude.
\end{proof}

\begin{lemma}\label{l2.2}
Let $\beta\in\C$ and $\Re\beta>-1/2$. Then
$\det T_{n}(v_\beta)\neq0$ for all $n\ge1$.
\end{lemma}
\begin{proof}
This follows from the formula
\bqn
\det T_n(v_\beta) &=& \frac{G(1+\beta)^2}{G(1+2\beta)}\cdot
\frac{G(1+n)\,G(1+2\beta+n)}{G(1+\beta+n)^2},
\nn
\eqn
which has been proved, e.g., in \cite{BS}.
\end{proof}

Let $\cP_{n,m}$ ($n\le m $) stand for the set of all trigonometric polynomials of the 
form
\bqn
p(t) &=& \sum_{k=n}^m p_k t^k.
\eqn
We also introduce the Hardy space
\bqn
\ovl{H^2(\T)}=\Big\{\;f\in L^2(\T)\;:\; f_k=0 \mbox{ for all } k>0\;\Big\}.
\eqn
which consists of those functions $f$ for which $\bar{f}\in H^2(\T)$.
Notice that $f\in H^2(\T)$ if and only if $\tilde{f}\in \ovl{H^2(\T)}$.
In the proof of the following we will use the identification (\ref{f.THid}).

\begin{proposition}\label{p2.3}
Let $\beta\in\C$ and $n\in\Z$. Then
\begin{itemize}
\item[(a)]
If $\Re\beta\in(-3/2,1/2]$, then 
\bqn
\dim\ker (I+H(u_{\beta-2n})) &=& \max\{0,-n\}.
\eqn
\item[(b)]
If $\Re\beta\in(-1/2,3/2]$, then 
\bqn
\dim\ker (I-H(u_{\beta-2n})) &=& \max\{0,-n\}.
\eqn
\end{itemize}
\end{proposition}
\begin{proof}
We will treat the operators $A=I+H(u_{\beta-2n})$ and $I-H(u_{\beta-2n})$
simultaneously. For sake of easy reference we will speak of Case A and Case B,
respectively.

Assume that $\Re\beta\in (-3/2,3/2)$ and let $f_+$ be in the kernel of 
$\ker(I\pm H(u_{\beta-2n}))$. Then
\bqn
f_+\mp t^{-2n-2}u_{\beta+1}\tilde{f}_+&=& t\iv f_-\in t\iv\ovl{H^2(\T)}.\nn
\eqn
We multiply with $\xi_{\beta+1}t^{n+1}$ and it follows that
\bqn
f_0:=\xi_{\beta+1}t^{n+1}f_+\mp t^{-n-1}\eta_{\beta+1}\tilde{f}_+&=& 
t^{n}\xi_{\beta+1} f_-\in t^n\ovl{H^1(\T)}.\nn
\eqn
Notice that $\xi_{\beta+1}\in\ovl{H^2(\T)}$. Obviously, $\tilde{f}_0=\mp f_0$.
Comparing the Fourier coefficients, it follows that the right hand side is zero
if $n<0$. 

We claim that the right hand side is also zero in the case $n\ge0$. In this case it follows 
first that $f_0=q_n$, where $q_n\in \cP_{-n,n}$ and $\tilde{q}_n=\mp q_n$. Hence
\bqn\label{eq:2-1}
t^{n+1}f_+\mp t^{-n-1}u_{\beta+1}\tilde{f}_+&=& \xi_{-\beta-1}q_n.
\eqn
We distinghish three cases.
\\
\underline{Case 1:} $\Re\beta\in (-3/2,-1/2)$. The last equation implies
\bqn\label{eq:2-2}
\eta_{-\beta-1}t^{n+1}f_+\mp \xi_{-\beta-1} t^{-n-1} \tilde{f}_+
&=&\xi_{-\beta-1}\eta_{-\beta-1} q_n,
\eqn
where $\xi_{-\beta-1}\in \ovl{H^2(\T)}$, $\eta_{-\beta-1}\in H^2(\T)$. For $k=-n,\dots,n$,
the $k$-th Fourier coefficient of $\xi_{-\beta-1}\eta_{-\beta-1}q_n$ is zero. 
This condition is equivalent to an equation 
$T_{2n+1}(\xi_{-\beta-1}\eta_{-\beta-1})\hat{q}_n=0$, where $\hat{q}_n$ is the vector consisting
of the Fourier coefficients of $q_n$. {}From Lemma \ref{l2.2} it follows that $q_n=0$.
\\
\underline{Case 2:} $\Re\beta\in [-1/2,1/2)$. Since $\xi_{-\beta-1}\not\in L^2(\T)$, equation 
(\ref{eq:2-1}) implies that $q_n(1)=0$. Write $q_n(t)=(1-t)q_{n-1}(t)$ with 
$q_{n-1}\in\cP_{-n,n-1}$. Multiplying (\ref{eq:2-2}) with $(1-t\iv)$, we obtain
\bqn\label{eq:2-3}
-\eta_{-\beta}t^{n}f_+\mp \xi_{-\beta} t^{-n-1} \tilde{f}_+
&=&\xi_{-\beta}\eta_{-\beta} q_{n-1}.
\eqn
Since $\xi_{-\beta}\in \ovl{H^2(\T)}$ and $\eta_{-\beta}\in H^2(\T)$,  for $k=-n,\dots,n-1$
the $k$-th Fourier coefficient of $\xi_{-\beta}\eta_{-\beta}q_{n-1}$ is zero. 
This condition lead us to 
$T_n(\xi_{-\beta}\eta_{-\beta})\hat{q}_{n-1}=0$ with $\hat{q}_{n-1}$ consisting of
the Fourier coefficients of $q_{n-1}$. Again Lemma \ref{l2.2} implies that $q_{n-1}=0$.
Hence $q_n=0$ as desired.
\\
\underline{Case 3:} $\Re\beta\in [1/2,3/2)$. Since $\xi_{-\beta}\not\in L^2(\T)$, equation 
(\ref{eq:2-1}) implies that $q_n(1)=q_n'(1)=0$. Write $q_{n}(t)=(1-t)(1-t\iv)q_{n-1}$
with $\cP_{-n+1,n-1}$. Multipliying (\ref{eq:2-2}) with $(1-t)(1-t\iv)$ it follows that
\bqn\label{eq:2-4}
-\eta_{-\beta+1}t^{n}f_+\pm \xi_{-\beta+1} t^{-n} \tilde{f}_+
&=&\xi_{-\beta+1}\eta_{-\beta+1} q_{n-1}.
\eqn
Notice that $\xi_{-\beta+1}\in\ovl{H^2(\T)}$ and $\eta_{-\beta+1}\in H^2(\T)$.
Similar as before, but now with the Toeplitz matrix $T_{2n-1}(\xi_{-\beta+1}\eta_{-\beta+1})$,
we obtain $q_n=0$.

After having proved that $q_n=0$ in all cases we can conclude that equations 
(\ref{eq:2-2}), (\ref{eq:2-3}) and (\ref{eq:2-4}) hold in all cases with the right hand side equal to zero.
It is now appropriate to distinguish again between several cases, but in a different way.
\\
\underline{Case (i):} $-1/2<\Re\beta <1/2$. From (\ref{eq:2-3}), i.e., 
$-\eta_{-\beta}t^{n}f_+\mp \xi_{-\beta} t^{-n-1} \tilde{f}_+=0$,
we obtain $f_+=0$ if $n\ge0$. If
$n<0$, then the general solution is $f_+=\eta_{\beta}p_n$ with $p_n\in\cP_{0,-2n-1}$
and $p_n(t)=\pm t^{-2n-1}\tilde{p}_n(t)$. Notice that $\eta_\beta\in H^2(\T)$ and that
the set of those polynomials $p_n$ is a linear space of dimension $-n$.
\\
\underline{Case (ii):} $-3/2<\Re\beta <-1/2$ and Case A. From (\ref{eq:2-2}), i.e., 
$\eta_{-\beta-1}t^{n+1}f_+ - \xi_{-\beta-1} t^{-n-1} \tilde{f}_+=0$, it follows $f_+=0$ if
$n\ge0$. If $n<0$, then the general solution is given by
$f_+=\eta_{\beta+1} p_n$ with $p_n\in \cP_{0,-2n-2}$ and $p_n(t)=t^{-2n-2} \tilde{p}_n(t)$.
The dimension of the space consisting of those polynomials $p_n$ is $-n$.
\\
\underline{Case (iii):} $1/2<\Re\beta <3/2$ and Case B. From (\ref{eq:2-4}), i.e., 
$\eta_{-\beta+1}t^{n}f_++ \xi_{-\beta+1} t^{-n} \tilde{f}_+=0$, we obtain $f_+=0$ in case
$n\ge0$. If $n<0$, then the general solution is given by
$f_+=\eta_{\beta-1} p_n$ with $p_n\in \cP_{0, -2n}$ and $p_n(t)= -t^{-2n}\tilde{p}_n(t)$.
The space of those polynomials is $-n$.
\\
\underline{Case (iv):} $\Re\beta=-1/2$. Here we proceed as in case (i) and obtain that
the solution is of the form $f_+=\eta_{\beta}p_n$ if $n<0$. However, $\eta_\beta\not\in L^2(\T)$, 
which implies that $p_n(t)=(1-t)p_{n-1}(t)$. Hence the general solution is
$f_+=\eta_{\beta+1}p_{n-1}$ with $p_{n-1}\in\cP_{0,-2n-2}$ and $p_{n-1}(t)=\mp t^{-2n-2} \tilde{p}_{n-1}(t)$.
The dimension of the space of all solutions is $-n$ in Case A and  $-n-1$ in Case B.
\\
\underline{Case (v):} $\Re\beta=1/2$. Here we proceed as in case (iii) and obtain that
the solution is of the form $f_+=\eta_{\beta-1}p_n$ if $n<0$. Since $\eta_{\beta-1}\not\in L^2(\T)$,
we can write $p_n(t)=(1-t)p_{n-1}(t)$. Hence the general solution is  
$f_+=\eta_{\beta}p_{n-1}$ with $p_{n-1}\in\cP_{0,-2n-1}$ and $p_{n-1}(t)=\pm t^{-2n-1} \tilde{p}_{n-1}(t)$.
The dimension of the space of all solutions is $-n$ both in Case A and in Case B.

Putting together the cases (i)-(v),  the statement of the proposition follows easily.
\end{proof}

\begin{theorem}\label{t2.4}
Let $\beta\in \C$. Then
\begin{itemize}
\item[(a)]
$I+H(u_{\beta})$ is invertible of $\el$ if and only if $\Re\beta<1/2$ and $\Re\beta\notin\frac{1}{2}+2\Z$,
\item[(b)]
$I-H(u_{\beta})$ is invertible of $\el$ if and only if $\Re\beta<3/2$ and $\Re\beta\notin\frac{3}{2}+2\Z$.
\end{itemize}
\end{theorem}
\begin{proof}
The Fredholm criteria is contained in Proposition \ref{p2.1}.
The dimension of the kernel and cokernel of $I\pm H(u_\beta)$ is given by Proposition
\ref{p2.3}. Notice that $H(u_\beta)^*=H(u_{\ovl{\beta}})$.
\end{proof}


\subsection{The determinants of the discrete operators}

For $\beta\in \C$ and $r\in[0,1)$ we introduce the functions 
\be
v_{\beta,r}(t):=(1-r/t)^{\beta}(1-rt)^{\beta},\qquad
u_{\beta,r}(t):=(1-r/t)^{-\beta} (1-rt)^{\beta},\qquad
t\in\T.
\ee
We will use these functions as approximations of the functions
$v_\beta$ and $u_\beta$. Recalling (\ref{f.uv-rel}) notice that we can write
$$
v_{\beta}(t)=(1-1/t)^{\beta}(1-t)^{\beta},\qquad
u_{\beta}(t)=(1-1/t)^{-\beta} (1-t)^{\beta},\qquad
t\in\T.
$$

\begin{lemma}\label{l2.5} \
\begin{itemize}
\item[(a)]
If $-3/2<\Re\beta<1/2$, then $(I+H(u_{\beta,r}))\iv\to (I+H(u_\beta))\iv$ strongly on $\el$ as
$r\to1$.
\item[(b)]
If $-1/2<\Re\beta<3/2$, then $(I-H(u_{\beta,r}))\iv\to (I-H(u_\beta))\iv$ strongly on $\el$ 
as $r\to1$.
\end{itemize}
\end{lemma}
\begin{proof}
We first remark that $u_{\beta,r}$ is bounded in the $L^\iy$-norm with respect to $r$
and converges in measure to $u_\beta$ as $r\to1$. Hence $H(u_{\beta,r})\to H(u_\beta)$
strongly on $\el$ as $r\to 1$. In order to prove the strong convergence of
the inverses of $I\pm H(u_{\beta,r})$ it is thus necessary and sufficient to show
that the sequence $I\pm H(u_{\beta,r})$ is stable.

Using the results of \cite{ES1} one can prove that 
$I\pm H(u_{\beta,r})$ is stable if and only if the operators
$$
I\pm H(u_{\beta})\quad\mbox{ and }\quad I\pm H(u_{-\beta,-1})
$$
are invertible, where $u_{-\beta,-1}(t):=u_{-\beta}(-t)$. 
Introducing the operator $W:\{x_n\}_{n=0}^\iy\in\el\mapsto\{(-1)^nx_n\}_{n=0}^\iy\in\el$ and 
noting that 
\be\label{f.Wflip}
W^2=I\quad\mbox{ and }\quad
WH(a)W=-H(b)
\ee
with $b(t):=a(-t)$, $t\in\T$, it follows that $I\pm H(u_{-\beta,-1})$
is invertible if and only if the operator $I\mp H(u_{-\beta})$ is invertible.
Applying Theorem \ref{t2.4} completes the proof.

In order to give some more details about the derivation of the above stability criterion from 
\cite{ES1} we rely on  the notation introduce there. What we have to do is to apply Theorem
2.1 and Theorem 2.2 of  \cite{ES1} in the setting where $k_\lambda$ equals the harmonic extension
(1.8).  The corresponding functions $K(x)$ and $f(e^{i\theta})$ (see (2.4)) evaluate to 
$K(x)=1/(\pi(1+x^2))$ and $f(e^{i\theta})=(\theta+\pi)/(2\pi)$, $|\theta|<\pi$. Hence the functions $a_\tau$
($\tau\in\T$) that are associated to each function $a\in PC$ (see (2.5)) are given by
$$
a_\tau(e^{i\theta})=a(\tau+0)\frac{\pi+\theta}{2\pi}+a(\tau-0)\frac{\pi-\theta}{2\pi},\qquad
-\pi<\theta<\pi.
$$
We need those functions in the case $a=\log u_\beta$.

The setting of Theorem 2.2 is with $A_\lambda=I\pm H(\exp(k_\lambda(\log u_\beta)))$
(since $k_\lambda (\log u_\beta)=\log u_{\beta,r}$ with $\lambda=-1/\log r$). The homomorphisms
evaluate with the help of Theorem 2.1 to
\begin{itemize}
\item[(i)] $\Psi_0[A_\lambda]=I\pm H(\exp(\log u_\beta))=I\pm H(u_\beta)$,
\item[(ii)] $\Psi_1[A_\lambda]=I\pm H(\exp(\log u_\beta)_1)=I\pm H(u_{-\beta,-1})$,
\item[] $\Psi_{-1}[A_\lambda]=I\mp H(\exp(\log u_\beta)_{-1})=I$,
\item[(iii)]  $\Psi_{\tau,\bar{\tau}}[A_\lambda]=
\twomat{I&0\\0&I}\pm \twomat{0& PM(\exp(\log u_\beta)_\tau)Q\\
QM(\exp\widetilde{(\log u_\beta)_\tau})P&0}=\twomat{I&0\\0&I}$.
\end{itemize}
Observe that $(\log u_\beta)_\tau$ is a constant except for $\tau\neq1$ since $u_\beta$
has only a discontinuity at $1$. For $\tau=1$ we have
$$
(\log u_\beta)_1(e^{i\theta})=-i\beta\pi \frac{\pi+\theta}{2\pi}+i\beta \pi\frac{\pi-\theta}{2\pi}
=-i\beta\theta, \qquad -\pi<\theta<\pi.
$$
Hence $\exp((\log u_\beta)_1(e^{i\theta}))=e^{-i\beta\theta}=u_{-\beta}(e^{i(\theta+\pi)})
=u_{-\beta,-1}(e^{i\theta})$, $|\theta|<\pi$, which settles (ii).

The invertibility of the operators in (i)--(iii) (which is necessary and sufficient for the stability of
the sequence $A_\lambda=I\pm H(u_{\beta,r})$) is nothing else than what has been stated above.
This completes the derivation.
\end{proof}

A simple conclusion of the previous lemma is the following result.

\begin{proposition}\label{p2.6}Let $n\ge1$ be fixed.
\begin{itemize}
\item[(a)]
If $-3/2<\Re\beta<1/2$, then
$$\det \Big[P_n (I+H(u_{\beta,r}))\iv P_n\Big]\to \det\Big[P_n(I+H(u_\beta))\iv P_n\Big],
\qquad r\to 1.$$
\item[(b)]
If $-1/2<\Re\beta<3/2$, then
$$\det \Big[P_n (I-H(u_{\beta,r}))\iv P_n\Big]\to \det \Big[P_n (I-H(u_\beta))\iv P_n\Big],
\qquad r\to 1.$$
\end{itemize}
\end{proposition}

Let $\W(\T)$ be the Banach algebra of all functions defined on the unit circle
which are of the form $a(t)=\sum_{n=-\iy}^\iy a_nt^n$ where
\bqn
\sum_{n=-\iy}^\iy |na_n|<\iy.
\eqn
Notice that $a\in \W(\T)$ implies that $H(a)$ is a trace class operator.
A Wiener-Hopf factorization in $\W(\T)$ is a representation of the form
\bqn
a(t)&=&a_-(t)a_+(t),\qquad t\in\T, 
\eqn
such that the functions $a_\pm$ and their inverses belong to $\W(\T)$ 
and such that the $n$-th Fourier coefficients of the functions
$a_+^{\pm1}(t)$ and $a_-^{\pm1}(t\iv)$ vanish for each $n<0$.
It is well-known \cite{BS} that a function $a\in\W(\T)$ possesses
a Wiener-Hopf factorization in $\W(\T)$ if and only if 
$a$ is nonzero on all of $\T$ and has winding number zero.
This is equivalent to the condition that $a$ possesses
a logarithm $\log a\in\W(\T)$.

Under this last condition we can define
\bqn
G[a] &:=&\exp\Big(\frac{1}{2\pi}\int_0^{2\pi} \log a(e^{i\theta})\,d\theta\Big)
\eqn
as the geometric mean of $a$.

\begin{proposition}\label{p2.5}
Let $a\in\W(\T)$ be an even function which possesses a Wiener-Hopf
factorization $a(t)=a_-(t)a_+(t)$. Define
$\psi(t)=\ta_+(t) a_+\iv(t)$.
Then $I\pm H(\psi)$ is invertible on $\el$ and
\bqn
\det \Big[ T_n(a)\pm H_n(a)\Big] &=& G[a]^n
\det \Big[P_n (I\pm H(\psi))\iv P_n\Big].
\eqn
\end{proposition}
\begin{proof}
First of all notice that from (\ref{f.Tab}) 
$$
\Big(T(a)\pm H(a)\Big)\Big(T(a\iv)\pm H(a\iv)\Big)=
\Big(T(a\iv)\pm H(a\iv)\Big)\Big(T(a)\pm H(a)\Big)=I.
$$
Moreover, using (\ref{f.Hab}) it follows that
$$
T(a\iv)\pm H(a\iv)= T(a_-\iv)(I\pm H(\psi))T(a_+\iv).
$$
Notice that also from  (\ref{f.Tab})
$$
T(a_\pm)T(a_\pm\iv)=T(a_\pm\iv)T(a_\pm)=I.
$$
Hence $I\pm H(\psi)$ is invertible and
$$
T(a)\pm H(a) = T(a_+)(I\pm H(\psi))\iv T(a_-).
$$
Now we multiply from the left and right with $P_n$, and observing
that $T(a_+)$ and $T(a_-)$ are lower and upper triangular
matrices we obtain
$$
T_n(a)\pm H_n(a) = T_n(a_+)\Big( P_n(I\pm H(\psi))\iv P_n\Big) T_n(a_-).
$$
Since $\det T_n(a_\pm) =([a_\pm]_0)^n$ and
$[a_+]_0[a_-]_0=\exp([\log a_+]_0+[\log a_-]_0)=G[a]$
we conclude the desired assertion.
\end{proof}

\begin{proposition}\label{p2.7} Let $n\ge1$ be fixed.
\begin{itemize}
\item[(a)]
If $-1/2<\Re\beta < 3/2$, then
\bqn\label{f.42}
D_n^+(\beta) &=&
\det \Big[P_n (I+H(u_{-\beta}))\iv P_n\Big].
\eqn
\item[(b)]
If $-3/2<\Re\beta < 1/2$, then
\bqn\label{f.43}
D_n^-(\beta) &=&
\det \Big[P_n (I-H(u_{-\beta}))\iv P_n\Big].
\eqn
\end{itemize}
\end{proposition}
\begin{proof}
We apply the previous proposition with $a(t)=v_{\beta,r}(t)$
(where $0\le r<1$). Noting that $a_+(t)=(1-rt)^\beta$,
$\psi(t)=u_{-\beta,r}(t)$ and $G[a]=1$
we obtain
$$
\det \Big[ T_n(v_{\beta,r})\pm H_n(v_{\beta,r})\Big] =
\det \Big[P_n (I\pm H(u_{-\beta,r}))\iv P_n\Big] 
$$
for all $\beta\in \C$.
For $\Re\beta>-1/2$, we have that $v_{\beta,r}\to v_\beta$ in the $L^1$-norm.
Hence the limit as $r\to1$ of the left hand side of the previous identity is (for $n$ fixed)
equal to $\det [T_n(v_\beta)\pm H_n(v_\beta)]$. From Proposition \ref{p2.6} 
we obtain the limit of the right hand side and thus the identities
(\ref{f.42}) and (\ref{f.43}).

In order to justify identity (\ref{f.43}) in the case $-3/2<\Re\beta\le-1/2$ in (b)
we argue by analyticity (see Theorem \ref{t.Ddiscrete}).
Notice that the analyticity of the determinant of the right hand side follows essentially
from the fact that the mapping $\beta\in\C\mapsto H(u_\beta)\in\cL(\el)$ is an
analytic operator-valued function.
\end{proof}

Obviously, the previous result in connection with Theorem \ref{t.Ddiscrete} 
allows us to determine the asymptotics of the determinants
$\det[P_n (I\pm H(u_{-\beta}))\iv P_n]$ as $n\to\iy$ in the case where
$-1/2<\pm\Re\beta<3/2$. This will be one of the cornerstones in the proof of the
main result (see Section \ref{s3.6}).


\subsection{The determinants of the continuous operators}
\label{sec:3:4}

In this subsection we will establish the continuous analogues
to the results of the previous subsection.
We introduce the functions
\be
\hv_{\beta,\eps}(x):=\left(\frac{x^2+\eps^2}{x^2+1}\right)^{\beta},\qquad
\hu_{\beta,\eps}(x):=\left(\frac{x-\eps i}{x-i}\right)^{-\beta}
\left(\frac{x+\eps i}{x+i}\right)^{\beta},
\ee
where $\beta\in\C$ and $\eps\in(0,1]$. These functions will
approximate the functions $\hv_\beta$ and $\hu_\beta$, respectively, which were defined
in (\ref{f.vuh}).

\begin{proposition}\label{p2.10} Let $R>0$ be fixed.
\begin{itemize}
\item[(a)]
If $-3/2<\Re\beta<1/2$, then
$$\det\Big[P_R(I+H(\hu_{\beta,\eps}))\iv P_R\Big]\to \det\Big[P_R(I+H(\hu_{\beta}))\iv P_R\Big],
\qquad \eps\to0.$$
\item[(b)]
If $-1/2<\Re\beta<3/2$, then
$$\det\Big[P_R(I-H(\hu_{\beta,\eps}))\iv P_R\Big]\to \det\Big[P_R(I-H(\hu_{\beta}))\iv P_R\Big],
\qquad \eps\to0.$$
\end{itemize}
\end{proposition}
\begin{proof}
We first write
$$
P_R(I\pm H(\hu_{\beta,\eps}))\iv P_R=
P_R\mp P_R(I\pm H(\hu_{\beta,\eps}))\iv H(\hu_{\beta,\eps}) P_R
$$
Noting that $W_R^2=P_R$, $W_RP_R=P_RW_R=W_R$, where $W_R:=H(e^{ixR})$, it follows that
$$
\det\Big[P_R(I\pm H(\hu_{\beta,\eps}))\iv P_R\Big]
=\det[P_R\mp A_\eps],
$$
where
$$
A_\eps:=W_R (I\pm H(\hu_{\beta,\eps}))\iv H(\hu_{\beta,\eps}) W_R.
$$
We claim that $H(\hu_{\beta,\eps}) W_R$ is a trace class operator, which converges
in the trace norm to $H(\hu_{\beta}) W_R$ as $\eps\to0$. 
In order to see this we apply the transform $\Phi$ (see (\ref{f.Phi})) and obtain
$$
\Phi[H(\hu_{\beta,\eps})W_R]=H(u_{\beta,r})H(h_R),
\qquad
\Phi[H(\hu_{\beta})W_R]=H(u_{\beta})H(h_R),
$$
where $h_R(t):=\exp(R(t-1)/(t+1))$, $r=(1-\eps)/(1+\eps)$. Let $f_1,f_2$ be smooth functions on $\T$
satisfying $f_1+f_2=1$ such that $f_1(t)$ vanishes on a neighborhood of $t=1$ and $f_2(t)$ 
vanishes on a neighborhood of $t=-1$. By applying (\ref{f.Hab}) we decompose
\bqn
H(u_{\beta,r}) H(h_R) &=& H(u_{\beta,r})\Big[T(f_1)+T(f_2)\Big]H(h_R)\nn\\
&=& \Big[H(u_{\beta,r} \tilde{f}_1)-T(u_{\beta,r})H(\tilde{f}_1)\Big]H(h_R)\nn\\
&&\mbox{}+H(u_{\beta,r})\Big[H(f_2 h_R)-H(f_2)T(\tilde{h}_R)\Big].\nn
\eqn
A similar identity where $u_{\beta,r}$ is replaced with $u_\beta$ can also be established.
It is easy to verify that $u_{\beta,r} \tilde{f}_1\to u_{\beta} \tilde{f}_1$
in the norm of $\W(\T)$ as $r\to1$ and that $f_2 h_R\in \W(\T)$. All Hankel operators
appearing within the brackets are trace class and $H(u_{\beta,r}\tilde{f}_1)\to
H(u_\beta\tilde{f}_1)$ in the trace norm.
Since $H(u_{\beta,r})\to H(u_\beta)$ and $T(u_{\beta,r})\to T(u_\beta)$
strongly on $\el$ (see Lemma \ref{l3.1}) we conclude that 
$H(u_{\beta,r})H(h_R)\to H(u_{\beta})H(h_R)$ in the trace norm. 
Hence $H(\hu_{\beta,\eps}) W_R$ 
is a trace class operator which converges in the trace norm to $H(\hu_{\beta}) W_R$.

Using Lemma \ref{l2.5} and the transform $\Phi$ we conclude that 
$(I\pm H(\hu_{\beta,\eps}))\iv$ converges strongly to
$(I\pm H(\hu_{\beta}))\iv$ as $\eps\to0$. Hence $A_\eps\to A$ in the trace norm where
$A$ is the trace class operator
$$ 
A:=W_R (I\pm H(\hu_{\beta}))\iv H(\hu_{\beta}) W_R.
$$
Thus we obtain $\det[P_R\mp A_\eps]\to \det[P_R\mp A]$,
which proves the assertion.
\end{proof}

Let $\W(\R)$ be the set of all functions $\ha$ defined on $\R$ such that 
$a\in\W(\T)$ where
\bqn\label{f.46}
a\left(\frac{1+ix}{1-ix}\right) &:=& \ha(x),\qquad x\in\R.
\eqn
Obviously, $\ha\in \W(\R)$ implies that $H(\ha)$ is a trace class operator on $L^2(\R_+)$
(see  (\ref{f.cd}) and (\ref{f.frac})).

For a function $\ha\in \W(\R)$ which possesses a logarithm 
$\log \ha\in L^1(\R)\cap \W(\R)$, the geometric means is well
defined by
\bqn\label{f.55}
G[\ha] &:=& \exp\left(\frac{1}{2\pi}\int_{-\iy}^\iy
\log \ha(x)\,dx\right).
\eqn
Notice that the logarithm is uniquely determined.

We say that $\ha(x)=\ha_-(x)\ha_+(x)$, $x\in\R$, is a Wiener-Hopf 
factorization in $\W(\R)$ if $a(t)=a_-(t)a_+(t)$, $t\in\T$, 
is a Wiener-Hopf factorization in $\W(\T)$, where the functions
$a$ and $a_\pm$ are defined according to (\ref{f.46}).

\begin{lemma}\label{l2.11a}
Let $a\in\W(\R)$ be a function which possesses a Wiener-Hopf
factorization $a(x)=a_-(x)a_+(x)$ in $\W(\R)$
and a logarithm $\log a\in L^1(\R)\cap \W(\R)$. 
Then $W_R(a_-)W_R(a_+)-P_R$ is a trace class 
operator on $L^2[0,R]$ and
\bqn
\det\Big[ W_R(a_-)W_R(a_+)\Big]=G[a]^R.
\eqn
\end{lemma}
\begin{proof}
We can assume without loss of generality that the factors are normalized
such that $(\log a_\pm)(\iy)=0$. Obviously, $\log a\in L^1(\R)\cap \W(\R)$
implies that $\log a\in L^2(\R)\cap \W(\R)$. Notice that 
$L^p(\R)\cap \W(\R)$ are Banach algebras without unit elements. 
Since the Riesz projection 
with respect to the upper half-plane is bounded on $L^2(\R)\cap \W(\R)$
it follows 
$$
\log a_\pm \in L^2(\R)\cap \W(\R),\qquad
a_\pm -1 \in L^2(\R)\cap \W(\R).
$$
Hence $W_R(a_\pm -1)$ are Hilbert-Schmidt
operators, while $W_R(a_\pm -1 -\log a_\pm)$ are trace class operators.
The latter is true since $a_\pm -1-\log a_\pm \in L^1(\R)$.
On the other hand $W_R(\log a)=W_R(\log a_+)+W_R(\log a_-)$
is also a trace class operator.

For Hilbert-Schmidt operators $K,L$ for which $K+L$ is a trace class operator
the identity
$$
\det \Big[(I+K)(I+L)\Big]=
\det \Big[(I+K)e^{-K}\Big]\det\Big[(I+L)e^{-L}\Big]\exp\Big[\trace(K+L)\Big]
$$
holds, which can be proved by an approximation argument.
We use this identity in the setting $K=W_R(a_-)-I$ and $L=W_R(a_+)-I$
and remark that the operator
\bqn
K+L &=& W_R(a_--1)+W_R(a_+-1)\nn\\
&=& W_R(a_--1-\log a_-)+W_R(a_+-1-\log a_+)+W_R(\log a)\nn
\eqn
is trace class. Noting that $W_R(a_\pm)=e^{W_R(\log a_\pm)}$, which implies
\bqn
\det\Big[(I+K)e^{-K}\Big] &=& \det W_R(a_-)e^{W_R(1-a_-)} =
\exp\Big[\trace\, W_R(\log a_-+1-a_-)\Big],
\nn\\
\det\Big[(I+L)e^{-L}\Big] &=& \det W_R(a_+)e^{W_R(1-a_+)} 
=\exp\Big[\trace\, W_R(\log a_++1-a_+)\Big],
\nn
\eqn
it follows that $\det \Big[(I+K)(I+L)\Big]$ equals the exponential of
$$
\trace \,W_R(\log a_-+1-a_-)+\trace\, W_R(\log a_++1-a_+)
+\trace(K+L)=\trace\, W_R(\log a).
$$
Since the trace of $W_R(\log a)$ is equal to $R$ times the Fourier transform of $\log a$
evaluated at the point $\xi=0$, the assertion follows easily.
\end{proof}

In regard to the proof of the previous lemma we remark that in general
$\log a_\pm \notin L^1(\R)$. In particular, $W_R(a_\pm)-P_R$ need not be trace class
operators.

\begin{proposition} \label{p2.11}
Let $a\in\W(\R)$ be an even function which possesses a Wiener-Hopf factorization
$a(x)=a_-(x)a_+(x)$ in $\W(\R)$. Suppose that $\log a\in L^1(\R)\cap\W(\R)$
and define $\psi(x)=\ta_+(x)a_+\iv(x)$. Then $I\pm H(\psi)$ is invertible
on $L^2(\R_+)$, and
\bqn
\det\Big[W_R(a)\pm H_R(a)\Big] &=& G[a]^R
\det\Big[P_R(I\pm H(\psi))\iv P_R\Big] .\nn
\eqn
\end{proposition}
\begin{proof}
We can prove in the same way as in Proposition \ref{p2.5} that 
$I\pm H(\psi)$ is invertible, and we derive the identity
$$
W(a)\pm H(a)=W(a_+)(I\pm H(\psi))\iv W(a_-).
$$
Since $a_\pm$ are appropriate Wiener-Hopf factors, we have 
$P_RW(a_+)=W_R(a_+)$ and $W(a_-)P_R=W_R(a_-)$. Hence
$$
W_R(a)\pm H_R(a)=W_R(a_+)\Big(P_R (I\pm H(\psi))\iv P_R\Big) W_R(a_-).
$$
The fact that $\psi\in \W(\R)$ implies that $P_R(I\pm H(\psi))\iv P_R$
is identity plus a trace class operator. Moreover, because the operators
$W_R(a_\pm)$ are invertible, we obtain 
\bqn
\det \Big[W_R(a)\pm H_R(a)\Big] &=&
\det \Big[W_R(a_-)W_R(a_+) P_R (I\pm H(\psi))\iv P_R\Big]
\nn\\[1ex]
&=&
\det\Big[ W_R(a_-)W_R(a_+)\Big]\cdot \det\Big[ P_R (I\pm H(\psi))\iv P_R\Big],\nn
\eqn 
which implies the assertion by employing Lemma \ref{l2.11a}.
\end{proof}

\begin{proposition}\label{p2.13}\
\begin{itemize}
\item[(a)] If $-1/2<\Re\beta<3/2$, then 
\bqn\label{f.49}
\hat{D}_R^+(\beta) &=&
e^{-\beta R}
\det\Big[P_R(I+H(\hu_{-\beta}))\iv P_R\Big].
\eqn
\item[(b)]
The function $\hat{D}_R^-(\beta)$ admits an analytic continuation onto
the set of all $\beta\in\C$ for which $\Re\beta>-3/2$.
Moreover, if $-3/2<\Re\beta<1/2$, then
\bqn\label{f.50}
\hat{D}_R^-(\beta) &=& e^{-\beta R}
\det\Big[P_R(I-H(\hu_{-\beta}))\iv P_R\Big].
\eqn
\end{itemize}
\end{proposition}
\begin{proof}
We already know that $\hat{D}_R^\pm(\beta)$ are analytic functions on the
set of all $\beta\in\C$ for which $\Re\beta>-1/2$. Moreover, the right hand side
in (\ref{f.50}) is analytic for $-3/2<\Re\beta<1/2$. This follows from the
fact that $H(\hat{u}_{\beta})$ is an operator-valued analytic function
in $\beta\in\C$ and that the inverses of $I-H(\hu_{-\beta})$
exist for $-3/2<\Re\beta<1/2$. 
Hence in order to prove statement (b) it suffices to prove the identity
(\ref{f.50}) for $-1/2<\Re\beta<1/2$.

We apply Proposition \ref{p2.11} with 
$a(x)=\hv_{\beta,\eps}(x)$. The corresponding Wiener-Hopf factors
are
$$
a_\pm(x)=\left(\frac{x\pm \eps i}{x\pm i}\right)^{\beta},
$$
whence we obtain $\psi(x) = \ta_+(x)a_+\iv(x)=\hu_{-\beta,\eps}(x)$.
Noting that $G[a]=e^{-\beta(1-\eps)}$ it follows that 
\bqn
\det\Big[ W_R(\hv_{\beta,\eps})\pm H_R(\hv_{\beta,\eps})\Big]
&=&  e^{-\beta R(1-\eps)}
\det\Big[P_R(I\pm H(\hat{u}_{\beta,\eps}))\iv P_R\Big].\nn
\eqn
Passing to the limit $\eps\to 0$ and applying Proposition \ref{p2.10}
we obtain the assertion.
\end{proof}

It is obvious from the previous proposition that we can determine the asymptotics
of $\hat{D}_R^\pm(\beta)$ from the asymptotics
of the determinant $\det\Big[ P_R(I\pm H(\hu_{-\beta}))\iv P_R\Big]$
and vice versa. This is the second ingredient in the proof
of the main result (see Section \ref{s3.6}).


\subsection{Asymptotic relation between discrete and continuous operators}
\label{s3.5}

In this section we are going to prove that (for certain fixed $\beta$)
$$
\det \Big[P_n(I\pm H(u_\beta))\iv P_n\Big]
\quad\sim\quad
\det \Big[P_R(I\pm H(\hu_\beta))\iv P_R\Big]
$$
as $n\to\iy$, $R\to\iy$ and $R=2n+O(1)$. 

We start with a couple of auxiliary results. The first result is one of the ingredients to the
proof of the Borodin-Okounkov identity as given in \cite{BW0,Bot}.

\begin{lemma}\label{l2.13}
Let $A$ be a trace class operator on a Hilbert space $H$ and assume that $I+A$ is invertible.
Let $P$ be a projection on $H$ and let $Q=I-P$. Then
\bqn
\det \Big[P(I+A)\iv P\Big] &=& \frac{\det(I+QAQ)}{\det(I+A)}.
\eqn
\end{lemma}
\begin{proof}
We write $(I+A)\iv=I-(I+A)\iv A$ and extend the operator appearing on the
left hand side in the operator determinant by the projection $Q$,
\bqn
P(I+A)\iv P +Q &=& I-P(I+A)\iv AP.\nn
\eqn
It follows that 
\bqn
\det \Big[P(I+A)\iv P\Big] &=& \det\Big[I-P(I+A)\iv A P\Big]\nn\\
&=& \det\Big[I-(I+A)\iv A P\Big]\nn\\
&=& \det(I+A)\iv\cdot\det\Big[I+A-AP\Big]\nn\\
&=& \det(I+A)\iv\cdot\det\Big[I+QAQ\Big],\nn
\eqn
which is the desired assertion.
\end{proof}

\begin{lemma}\label{l2.14}
For $-1<\sigma<1$, the trace norm of the integral operator with the kernel 
$$
k(x,y)=\frac{f_1(x)f_2(y)}{x+y}
$$
on $L^2(M)$, where $M\subset \R_+$, is at most a constant times
the square root of 
$$
\int_{M}|f_1(x)|^2\,\frac{dx}{x^{1+\sigma}}\cdot \int_{M}|f_2(x)|^2\,\frac{dx}{x^{1-\sigma}}.
$$
\end{lemma}
\begin{proof}
We can write this operator as a product $K_1K_2$ where
$K_1:L^2(\R_+)\to L^2(M)$ and $K_2:L^2(M)\to L^2(\R_+)$
have the kernels
$$
k_1(x,\eta)=f_1(x)e^{-x\eta}\eta^{\sigma/2},\qquad
k_2(\xi,y)=f_2(y)e^{-y\xi}\xi^{-\sigma/2}.
$$
The operators $K_1$ and $K_2$ are Hilbert-Schmidt and their norms can 
be estimated appropriately.
\end{proof}

Let $K_{\beta,\eps,n}$ and $\hK_{\beta,\eps,R}$ be the integral operators
on $L^2[\eps,1]$ with the kernels
\bqn
K_{\beta,\eps,n}(x,y) &=&
-\frac{\sin(\pi\beta)}{\pi}
\left(\frac{(1+x)(x-\eps)}{(1-x)(x+\eps)}
\frac{(1+y)(y-\eps)}{(1-y)(y+\eps)}\right)^{\beta/2} 
\left(\frac{1-x}{1+x} \right)^{2n}\frac{1}{x+y},
\\
\hK_{\beta,\eps,R}(x,y) &=&
-\frac{\sin(\pi\beta)}{\pi}
\left(\frac{(1+x)(x-\eps)}{(1-x)(x+\eps)}
\frac{(1+y)(y-\eps)}{(1-y)(y+\eps)}\right)^{\beta/2} 
\frac{e^{-2Rx}}{x+y}.
\eqn

\begin{proposition}\label{p2.15}
Let $-1<\Re\beta<1$. Then $K_{\beta,\eps,n}$ and $\hK_{\beta,\eps,R}$ are
trace class operators on $L^2[\eps,1]$ and 
\bqn
\det(I\pm Q_nH(u_{\beta,r})Q_n) &=& \det(I \pm K_{\beta,\eps,n}),
\quad \textstyle r=\frac{1-\eps}{1+\eps},\label{f2.60}\\
\det(I\pm Q_RH(\hu_{\beta,\eps})Q_R) &=& \det(I \pm \hK_{\beta,\eps,R}),
\label{f2.61}
\eqn
where $Q_n=I-P_n$ and $Q_R=I-P_R$.
\end{proposition}
\begin{proof}
The fact that $K_{\beta,\eps,n}$ and $\hK_{\beta,\eps,R}$ are
trace class operators follows from Lemma \ref{l2.14}
(with $\sigma=0$).

Let us first prove identity (\ref{f2.60}).
The operator $Q_nH(u_{\beta,r})Q_n$ can be identified with the matrix kernel
\bqn
k(j,k) &=&
\frac{1}{2\pi i} \int_{\T} \left(\frac{1-rt}{1-rt\iv}\right)^\beta t^{-2-j-k-2n}\,dt\nn\\
&=&
\frac{1}{2\pi i} \int_{\T} \left(\frac{1-rt\iv}{1-rt}\right)^\beta t^{j+k+2n}\,dt.
\hspace*{10ex}\nn\\
&=&
-\frac{1}{\pi}\int_{-\iy}^\iy
\left(\frac{(\xi+i)(\xi-i\eps)}{(\xi-i)(\xi+i\eps)}\right)^\beta
\left(\frac{i-\xi}{i+\xi}\right)^{j+k+2n}\frac{d\xi}{(i+\xi)^2}.\nn
\eqn
Therein we have employed first the substitution $t\mapsto t\iv$ and then 
$t=\frac{i-\xi}{i+\xi}$, $r=\frac{1-\eps}{1+\eps}$.
The integrand is analytic in the upper
half-plane cut along the segment $[i\eps,i]$.
We deform the path of integration to this segment described back and forth. The expression
in parentheses is real and negative. 
The limit of its argument from the left equals $-\pi$ and from right equals $\pi$.
We obtain (with the substitution $\xi=i\eta$)
\bqn
k(j,k) &=&
-\frac{2\sin(\pi \beta)}{\pi} \int_{\eps}^1 
\left(\frac{(1+\eta)(\eta-\eps)}{(1-\eta)(\eta+\eps)}\right)^\beta 
\left(\frac{1-\eta}{1+\eta}\right)^{j+k+2n}\,\frac{d\eta}{(1+\eta)^2}.
\eqn
This operator can be written as $UV$ where
$U:L^2[\eps,1]\to\ell^2(\Z_+)$ and $V:\ell^2(\Z_+)\to L^2[\eps,1]$ are given by
\bqn
U(j,\xi) &=& -\frac{2\sin(\pi \beta)}{\pi} \left(\frac{1-\xi}{1+\xi}\right)^{j-\beta/2}
\left(\frac{\xi-\eps}{\xi+\eps}\right)^{\beta/2}
\frac{1}{(1+\xi)},\nn\\
V(\eta,k) &=& \left(\frac{1-\eta}{1+\eta}\right)^{2n+k-\beta/2}
\left(\frac{\eta-\eps}{\eta+\eps}\right)^{\beta/2} \frac{1}{(1+\eta)}.
\eqn
Under the assumption $-1<\Re\beta<1$, the operators $U$ and $V$ are Hilbert-Schmidt. 
The operator $VU$ is the integral operator with the kernel 
\bqn
h(\eta,\xi) &=&
-\frac{2\sin(\pi \beta)}{\pi} 
\left(\frac{(\eta-\eps)(\xi-\eps)}{(\eta+\eps)(\xi+\eps)}\right)^{\beta/2}
\frac{1}{(1+\xi)(1+\eta)}
\sum_{k=0}^\iy
\left(\frac{1-\eta}{1+\eta}\right)^{k+2n-\beta/2}
 \left(\frac{1-\xi}{1+\xi}\right)^{k-\beta/2}
\nn\\
&=&
-\frac{\sin(\pi \beta)}{\pi} 
\left(\frac{(\eta-\eps)(\xi-\eps)}{(\eta+\eps)(\xi+\eps)}
\frac{(1-\eta)(1-\xi)}{(1+\eta)(1+\xi)}\right)^{\beta/2}
\left(\frac{1-\eta}{1+\eta}\right)^{2n}
\frac{1}{\xi+\eta}.
\nn
\eqn
Hence $VU=K_{\beta,\eps,n}$.

Now we turn to the proof of (\ref{f2.61}).
Since $\hu_{\beta,\eps}-1\in L^2(\R)$, the operator $Q_R H(\hu_{\beta,\eps})Q_R$ 
can be identified with an integral operator with the kernel
\bqn
k(x,y)&=&\lim_{M\to\iy}\frac{1}{2\pi}\int_{-M}^{M} 
\left[\left(\frac{(\xi-i)(\xi+\eps i)}{(\xi+i)(\xi-\eps i)}\right)^\beta-1\right]
e^{-i\xi(2R+x+y)}\,d\xi\nn\\
&=& \lim_{M\to\iy}\frac{1}{2\pi}\int_{-M}^{M} 
\left[\left(\frac{(\xi+i)(\xi-\eps i)}{(\xi-i)(\xi+\eps i)}\right)^\beta -1\right]
e^{i\xi(2R+x+y)}\,d\xi.\nn
\eqn
The integrand is analytic in the upper half-plane cut along the segment $[i\eps,i]$ and decays
as $O(\xi\iv e^{-2R\,\Im \xi})$ as $\xi\to\iy$, $\Im\xi\ge0$. We deform the path of integration 
to this segment described back and forth. The expression in parentheses is real and negative. 
The limit of its argument from the left equals $-\pi$ and from right equals $\pi$. Hence we obtain
\bqn
k(x,y)&=& -\frac{\sin(\pi\beta)}{\pi}
\int_{\eps}^1 \left(\frac{(1+\eta)(\eta-\eps)}{(1-\eta)(\eta+\eps)}\right)^\beta 
e^{-(2R+x+y)\eta}\,d\eta.\nn
\eqn
This operator can be written as a product $UV$, where
$U:L^2[\eps,1]\to L^2(\R_+)$ and $V:L^2(\R_+)\to L^2[\eps,1]$ are given by
\bqn
U(x,\xi) &=& -\frac{\sin(\pi\beta)}{\pi} 
\left(\frac{(1+\xi)(\xi-\eps)}{(1-\xi)(\xi+\eps)}\right)^{\beta/2}e^{-x\xi},
\nn\\
V(\eta,y) &=& \left(\frac{(1+\eta)(\eta-\eps)}{(1-\eta)(\eta+\eps)}\right)^{\beta/2} 
e^{-(2R+y)\eta}.\nn
\eqn
Under the assumption $-1<\Re\beta<1$, the operators $U$ and $V$ are 
Hilbert-Schmidt operators. The operator $VU$ has the kernel
\bqn
h(\eta,\xi) &=&  -\frac{\sin(\pi\beta)}{\pi} 
\left(\frac{(1+\eta)(\eta-\eps)}{(1-\eta)(\eta+\eps)}
\frac{(1+\xi)(\xi-\eps)}{(1-\xi)(\xi+\eps)}\right)^{\beta/2}
\int_{0}^\iy
e^{-(2R+x)\eta-x\xi}\,dx\nn\\
&=&
-\frac{\sin(\pi\beta)}{\pi}
\left(\frac{(1+\eta)(\eta-\eps)}{(1-\eta)(\eta+\eps)}
\frac{(1+\xi)(\xi-\eps)}{(1-\xi)(\xi+\eps)}\right)^{\beta/2}
\frac{e^{-2R\eta}}{\eta+\xi},\nn
\eqn
which is the operator $\hK_{\beta,\eps,R}$.
\end{proof}

\begin{proposition}\label{p2.16}
Let $-1<\pm \Re\beta<1/2$. Then 
\bqn
\frac{\det P_R(I\pm H(\hat{u}_\beta))\iv P_R}
{\det P_n(I\pm H(u_\beta))\iv P_n}
&=&
\lim\limits_{\eps\to 0}
\frac{\det (I\pm \hat{K}_{\beta,\eps,R})}
{\det (I\pm K_{\beta,\eps,n})}.\label{f3.P}
\eqn
\end{proposition}
\begin{proof}
Applying Lemma \ref{l2.13} with $P=P_n$, $A=\pm H(u_{\beta,r})$, and
$P=P_R$, $A=\pm H(\hu_{\beta,\eps})$, respectively,
and Proposition \ref{p2.15}, it follows that
\bqn
\det  \Big[P_n(I\pm H(u_{\beta,r}))\iv P_n\Big] &=&
\frac{\det (I\pm K_{\beta,\eps,n})}{\det(I\pm H(u_{\beta,r}))},
\label{f3.Pn}
\\
\det  \Big[P_R(I\pm H(\hat{u}_{\beta,\eps}))\iv P_R\Big] &=&
\frac{\det (I\pm \hK_{\beta,\eps,R})}{\det(I\pm H(\hat{u}_{\beta,\eps}))},
\label{f3.PR}
\eqn
where $r=\frac{1-\eps}{1+\eps}$. By (\ref{f.cd}) and (\ref{f.frac}) the operators 
$H(u_{\beta,r})$ and $H(\hu_{\beta,\eps})$ are unitarily equivalent.
The invertibility of $I\pm H(u_{\beta,r})$ for $r$ sufficiently close to $1$ follows
from Lemma \ref{l2.5}. Hence the fractions on the right hand side of
(\ref{f3.Pn}) and (\ref{f3.PR}) are well defined for $r\to 1$ and $\eps\to 0$.

In fact,  one can even say more. From Proposition \ref{p2.5} with $\psi$ chosen
as in the proof Proposition \ref{p2.7}, it follows that 
$I\pm H(u_{\beta,r})$ is invertible for all $r\in[0,1)$.
Similarly, from Proposition \ref{p2.11} with $\psi$ chosen as in in proof of
Proposition \ref{p2.13}, it follows that 
$I\pm H(\hat{u}_{\beta,\eps})$ is invertible for all $\eps>0$.

Taking the quotient of (\ref{f3.Pn}) and (\ref{f3.PR}) and passing 
to the limit $\eps\to0$ we obtain the  desired assertion by using 
Proposition  \ref{p2.6} and \ref{p2.10}.
\end{proof}

One remark is in order concerning the non-vanishing of the denominators
of the fractions in (\ref{f3.P}). First of all, a careful examination of the expression for
$D_n^\pm(\beta)$ as stated in the proof of Theorem \ref{t.Ddiscrete} combined with the
exact formulas of Proposition \ref{p2.7} imply that the determinants
$\det P_n(I\pm H(u_\beta))\iv P_n$ are nonzero for all $n\ge 1$ and $\beta$ satisfying
$-3/2<\pm\Re\beta<1/2$. From Proposition \ref{p2.6} we can conclude that 
$\det P_n(I\pm H(u_{\beta,r}))\iv P_n$ are nonzero for $r$ sufficiently close to $1$.
Formula (\ref{f3.Pn}) now implies that also $\det(I\pm K_{\beta,\eps,n})$ is nonzero
for $\eps\to0$.

Our next step is to determine the limit $\eps\to 0$ on the
right hand sides of (\ref{f3.P}).
Before we are able to do this, we have to establish a couple
of auxiliary results. Some of them will be needed only later on
in order to analyze the expression which is obtained for 
the limit.

Let $K_\beta^0$, $K_{\beta,n}$ and $\hat{K}_{\beta,R}$ stand for the integral
operators on $L^2[0,1]$ with the following kernels:
\bqn
K^0_\beta(x,y) &=& -\frac{\sin(\pi\beta)}{\pi}\frac{1}{x+y},
\\
K_{\beta,n}(x,y) &=& -\frac{\sin(\pi\beta)}{\pi}
\left(\frac{1-x}{1+x}\right)^{2n-\beta/2}
\left(\frac{1-y}{1+y}\right)^{-\beta/2}\frac{1}{x+y},
\\
\hat{K}_{\beta,R}(x,y) &=& -\frac{\sin(\pi\beta)}{\pi}
\left(\frac{1-x}{1+x}\right)^{-\beta/2}
\left(\frac{1-y}{1+y}\right)^{-\beta/2}\frac{e^{-2Rx}}{x+y}.
\eqn
Moreover, let $H_\beta^0$ and $H_\beta$ stand for the integral operators with 
the following kernels on $L^2[1,\iy)$,
\bqn
H^0_\beta(x,y) &=& -\frac{\sin(\pi\beta)}{\pi}\frac{1}{x+y},
\\
H_{\beta}(x,y) &=& -\frac{\sin(\pi\beta)}{\pi}
\left(\frac{x-1}{x+1}\right)^{\beta/2}\left(\frac{y-1}{y+1}\right)^{\beta/2}
\frac{1}{x+y}.
\eqn
Finally, let $Y_\eps$ stand for the unitary operator 
$$
f(x)\in L^2[\eps,1]\mapsto \sqrt{\eps}f(\eps x)\in L^2[1,\eps\iv],
$$
and let $\Pi_{[a,b]}$ stand for the projections operator 
$f(x)\mapsto \chi_{[a,b]}(x)f(x)$, which is thought of acting
on appropriate spaces $L^2(M)$, $M\subset\R$.
 
In what follows we will prove that the above integral operators are bounded 
and that certain differences between them are even trace class. Moreover, certain
invertibility results will be established, too.

\begin{lemma}\label{l3.19}
The operators $K_\beta^0$ and $H_\beta^0$ are bounded. Moreover,
for $\Re\beta\notin \pm 1/2+2\Z$, the operators 
$I\pm K_\beta^0$ and $I\pm H_\beta^0$ are invertible,
$$
\Big(I\pm \Pe K_\beta^0 \Pe\Big)\iv\to 
(I\pm K_\beta^0)\iv,\qquad \eps\to0,
$$
strongly on $L^2[0,1]$, and
$$
\Big(I\pm Y_\eps \Pe K_\beta^0 \Pe Y_\eps^*\Big)\iv\to
(I\pm H^0_\beta)\iv,\qquad \eps\to0,
$$
strongly on $L^2[1,\iy)$.
\end{lemma}
\begin{proof}
The operator on $L^2(\R_+)$ with the kernel $\pi\iv(x+y)\iv$ is 
a bounded selfadjoint operator with spectrum equal to $[0,1]$. Indeed,
by a substitution $x\mapsto e^{-x}$, $y\mapsto e^{-y}$ it is easily seen that 
this operator is unitary equivalent to the integral operator on $L^2(\R)$ 
with the kernel
$(2\pi)\iv\sech((x-y)/2)$. This is a convolution operator 
with the symbol $\sech(\pi\xi)$, $\xi\in\R$, and thus its spectrum equals
$[0,1]$.

The restriction of this operator onto the spaces $L^2[0,1]$,
$L^2[\eps,1]$, $L^2[1,\eps\iv]$ and $L^2[1,\iy)$ are also bounded selfadjoint
operators with the spectrum contain in (in fact, equal to) the interval
$[0,1]$.
 
Hence under the above conditions on the parameter $\beta$, the operators 
$$
I\pm K^0_\beta,\quad
I\pm H^0_\beta,\quad
I\pm \Pe K^0_\beta \Pe,\quad
I\pm Y_\eps \Pe K^0_\beta \Pe Y_\eps^*
$$
are all bounded, and (in the last two cases) the norms of their inverses do not depend on $\eps$.
It remains to observe that 
$$
\Pe K_\beta^0 \Pe\to K_\beta^0 
\quad\mbox{and}\quad
Y_\eps \Pe K_\beta^0 \Pe Y_\eps^*=
\Pei H_\beta^0\Pei \to
H^0_\beta
$$ 
strongly on $L^2[0,1]$ and $L^2[1,\iy)$, respectively, as $\eps\to0$.
\end{proof}

The following lemma shows, in particular, that the operators 
$K_{\beta,n}$, $\hat{K}_{\beta,R}$ and $H_\beta$ are bounded for certain $\beta$.

\begin{lemma}\label{l3.20}
If $\Re\beta<1$, then the operators
$$
K_{\beta,n}-K_\beta^0\quad\mbox{ and }\quad
\hK_{\beta,R}-K_\beta^0
$$
are trace class operators and 
$$
K_{\beta,n}-\hK_{\beta,R}\to0
$$
in the trace norm as $R\to\iy$, $n\to\iy$, $R=2n+O(1)$.
If $\Re\beta>-1$, then 
$$H_\beta-H_\beta^0
$$ is a trace class operator.
\end{lemma}
\begin{proof}
The assertion that $K_{\beta,n}-K_\beta^0$ and $K_{\beta,R}-K_\beta^0$
are trace class operators can be proved by considering 
``intermediate'' operators with the kernel
$$
-\frac{\sin(\pi\beta)}{\pi}
\left(\frac{1-y}{1+y}\right)^{-\beta/2}\frac{1}{x+y},
$$
and by applying Lemma \ref{l3.19}. Similarly, the fact that 
$H_\beta-H_\beta^0$ is trace class can be proved by introducing
the operator with the kernel
$$
-\frac{\sin(\pi\beta)}{\pi}
\left(\frac{1-y}{1+y}\right)^{\beta/2}\frac{1}{x+y}.
$$

Finally, the trace norm of $K_{\beta,n}-\hK_{\beta,R}$ can be
estimated by a constant times the square-root of
$$
\int_0^1\left(\frac{1-y}{1+y}\right)^{-\Re\beta}\,\frac{dy}{y^{1/2}}\cdot
\int_0^1\left(\frac{1-x}{1+x}\right)^{-\Re\beta}\left|
\left(\frac{1-x}{1+x}\right)^{2n}-e^{-2Rx}\right|^2\frac{dx}{x^{3/2}}.
$$
The first integral is finite, and the second one can be split (for each $0<\delta<1$)
into an integral from $0$ to $\delta$ and an integral from $\delta$ to $1$.
The integral from $\delta$ to $1$ is finite and converges to zero as
$n,R\to\iy$. In the integral from $0$ to $\delta$ we estimate the first
term in the integrand by a constant (depending on $\delta$) and make a substitution 
$x\mapsto x/(4n)$ to obtain an upper estimate
$$
C_{\delta} n^{1/2}\int_0^{4n\delta}
\left|
\left(\frac{1-x/(4n)}{1+x/(4n)}\right)^{2n}-e^{-xR/(2n)}\right|^2\frac{dx}{x^{3/2}}.
$$
This equals
$$
C_\delta n^{1/2}\int_0^{4n\delta}
\left|
e^{-x+O(x^2/n)}-e^{-x+O(x/n)}\right|^2\frac{dx}{x^{3/2}}.
$$
Omitting the constant $C_\delta$, we split this integral into 
$$
n^{1/2}\int_0^{n^{1/3}}
\left|
e^{-x+O(x^2/n)}-e^{-x+O(x/n)}\right|^2\frac{dx}{x^{3/2}}
=
n^{1/2}\int_0^{n^{1/3}}e^{-2x}
O\left(\frac{x^2}{n^{4/3}}\right)\frac{dx}{x^{3/2}}
= O(n^{-5/6})
$$
and
$$
 n^{1/2}\int_{n^{1/3}}^{4n\delta}
\left|
e^{-x+O(x^2/n)}-e^{-x+O(x/n)}\right|^2\frac{dx}{x^{3/2}}
=
n^{1/2}\int_{n^{1/3}}^{4n\delta}
\left|
e^{-x+O(\delta x)}-e^{-x+O(\delta)}\right|^2\frac{dx}{x^{3/2}}
=
 O(e^{-n^{1/3}}),
$$
where the last estimate holds under the assumption that $\delta$ is
chosen small enough to guarantee that $O(x\delta)\le x/2$. Collecting all
terms, this proves the convergence of $K_{\beta,n}-\hK_{\beta,R}$
in the trace norm.
\end{proof}

\begin{lemma}\label{l3.21}
If $-3/2<\Re\beta<1/2$, then the inverses of 
$$
I+K_{\beta,n}\quad\mbox{ and }\quad
I+\hK_{\beta,R}
$$
exist for sufficiently large $n$ and $R$, respectively, and are uniformly bounded.
If $-1/2<\Re\beta<1$, then the inverses of 
$$
I-K_{\beta,n}\quad\mbox{ and }\quad
I-\hK_{\beta,R}
$$
exist for sufficiently large $n$ and $R$, respectively, and are uniformly bounded.
\end{lemma}
\begin{proof}
We prove the statements only for the case of the operators
$\hK_{\beta,R}$. The proof in the case of $K_{\beta,n}$ is analogous.
Introduce the operator $\hK'_{\beta,R}$ with the kernel
$$
\hK'_{\beta,R}(x,y) = -\frac{\sin(\pi \beta)}{\pi}\frac{e^{-2Rx}}{x+y}.
$$
The difference $\hK_{\beta,R}-\hK'_{\beta,R}$ can be estimated in the trace norm
by a constant times the square-root of the integrals 
$$
\int_0^1\left|\left(\frac{1-x}{1+x}\right)^{-\beta/2}-1\right|^2e^{-4Rx}\,\frac{dx}{x^{3/2}}\cdot
\int_0^1\left(\frac{1-y}{1+y}\right)^{-\Re\beta}\frac{dy}{y^{1/2}}.
$$
$$+
\int_0^1e^{-4Rx}\,\frac{dx}{x^{1/2}}\cdot
\int_0^1\left|\left(\frac{1-y}{1+y}\right)^{-\beta/2}-1\right|^2\frac{dy}{y^{3/2}}
$$
by using Lemma \ref{l2.14}. These terms converge to zero as $R\to\iy$ (under the assumption $\Re\beta<1$). 
Thus it is sufficient to prove that the inverses of $I\pm \hK'_{\beta,R}$ are
uniformly bounded. Now notice that $\hK'_{\beta,R}=A^2_R K_{\beta}^0$, where $A_R$ is the
multiplication operator with the symbol $e^{-Rx}$. Since $A_R$ is uniformly bounded,
the well-known relationship between the inverses of $I\pm AB$ and $I\pm BA$
implies that the remaining problem is reduced to showing that
the inverses of $I\pm A_R K_\beta^0 A_R$ are uniformly bounded. 
It remains to observe that 
$A_R=A_R^*$, $A_RA_R^*\le I$ and that the operator with the kernel 
$1/(\pi(x+y))$ (i.e., $K_\beta^0$ without the sine-factor)
is selfadjoint with its spectrum contained in $[0,1]$.
The proof can now be completed as in Lemma \ref{l3.19}.
\end{proof}

\begin{lemma}\label{l3.22}
Let $-1<\Re\beta<1$. Then
\bqn
\Pr K_{\beta,\eps,n} &=& \Pe  K_{\beta,n}\Pe + o_1(1),
\label{f.as1}\\[1ex]
\Pr\hat{K}_{\beta,\eps,R} &=&\Pe \hat{K}_{\beta,R}\Pe+o_1(1),
\label{f.as2}\\[1ex]
Y_\eps \Per K_{\beta,\eps,n}Y_\eps^*
&=& \Pei H_{\beta}\Pei+o_1(1),
\label{f.as3}\\[1ex]
Y_\eps \Per \hat{K}_{\beta,\eps,R}Y_\eps^*
&=& \Pei H_{\beta}\Pei+o_1(1)
\label{f.as4}
\eqn
as $\eps\to0$, where $o_1(1)$ stands for a sequence of operator converging to 
zero in the trace norm.
\end{lemma}
\begin{proof}
We are going to prove only the identities involving $\hK_{\beta,\eps,R}$.
The assertions involving $K_{\beta,\eps,n}$ can be proved analogously.

As to identity (\ref{f.as2}), we have to show that the integral operator on $L^2[\eps,1]$
with the kernel
$$
\chi_{[\sqrt{\eps},1]}
\left(\frac{(1-x)(1-y)}{(1+x)(1+y)}\right)^{-\beta/2}
\left[\left(\frac{(x-\eps)(y-\eps)}{(x+\eps)(y+\eps)}\right)^{\beta/2}-1
\right]\frac{e^{-2xR}}{x+y}
$$
converges in  the trace norm to zero.
We split this kernel into the sum of the kernels
$$
\chi_{[\sqrt{\eps},1]}
\left(\frac{(1-x)(1-y)}{(1+x)(1+y)}\right)^{-\beta/2}
\left(\frac{x-\eps}{x+\eps}\right)^{\beta/2}
\left[\left(\frac{y-\eps}{y+\eps}\right)^{\beta/2}-1
\right]\frac{e^{-2xR}}{x+y}
$$
and
$$
\chi_{[\sqrt{\eps},1]}
\left(\frac{(1-x)(1-y)}{(1+x)(1+y)}\right)^{-\beta/2}
\left[\left(\frac{x-\eps}{x+\eps}\right)^{\beta/2}-1
\right]\frac{e^{-2xR}}{x+y}.
$$
The first of these kernels can be estimated by
$$
\int_{\sqrt{\eps}}^1
\left(\frac{1-x}{1+x}\right)^{-\Re\beta}
\left(\frac{x-\eps}{x+\eps}\right)^{\Re\beta}
\,\frac{dx}{x^{3/2}}\cdot
\int_\eps^1
\left(\frac{1-y}{1+y}\right)^{-\Re\beta}
\left|\left(\frac{y-\eps}{y+\eps}\right)^{\beta/2}-1
\right|^2\,\frac{dy}{y^{1/2}},
$$
and the second one can be estimated by
$$
\int_{\sqrt{\eps}}^1
\left(\frac{1-x}{1+x}\right)^{-\Re\beta}
\left|\left(\frac{x-\eps}{x+\eps}\right)^{\beta/2}-1\right|^2
\,\frac{dx}{x^{3/2}}\cdot
\int_\eps^1
\left(\frac{1-y}{1+y}\right)^{-\Re\beta}\,\frac{dy}{y^{1/2}}.
$$
We split off from all these integrals integrals from $1/2$
to $1$ in order to get rid off the sigularity at $1$.
In the remaining integrals (from $\sqrt{\eps}$ to $1/2$
and $\eps$ to $1/2$, resp.) we make a substitution
$x=\sqrt{\eps}z$ and $y=\eps z$, respectively.
Collecting all terms we obtain
$(O(1)+O(\eps^{-1/4}))(O(\eps^2)+O(\eps^{1/2}))=O(\eps^{1/4})$
for the first expression and
$(O(\eps)+O(\eps^{3/4}))O(1)=O(\eps^{3/4})$
for the second expression. Hence both terms 
converge to zero as $\eps\to 0$.

As to (\ref{f.as4}), we have to prove that 
the integral operator on $L^2[1,\eps\iv]$ with the
kernel
$$
\chi_{[1,1/\sqrt{\eps}]}(x)
\left(\frac{(x-1)(y-1)}{(x+1)(y+1)}\right)^{\beta/2}
\left[\left(\frac{(1-x\eps)(1-y\eps)}{(1+x\eps)(1+y\eps)}\right)^{-\beta/2}
e^{-2Rx\eps}-1\right]\frac{1}{x+y}
$$
tends to zero in the trace norm. We split this kernel into
$$
\chi_{[1,1/\sqrt{\eps}]}(x)
\left(\frac{(x-1)(y-1)}{(x+1)(y+1)}\right)^{\beta/2}
\left(\frac{1-x\eps}{1+x\eps}\right)^{-\beta/2}
e^{-2Rx\eps}
\left[
\left(\frac{1-y\eps}{1+y\eps}\right)^{-\beta/2}
-1\right]\frac{1}{x+y}
$$
and
$$
\chi_{[1,1/\sqrt{\eps}]}(x)
\left(\frac{(x-1)(y-1)}{(x+1)(y+1)}\right)^{\beta/2}
\left[\left(\frac{1-x\eps}{1+x\eps}\right)^{-\beta/2}
e^{-2Rx\eps}-1\right]\frac{1}{x+y}.
$$
These kernels can be estimated by
$$
\int_{1}^{1/\sqrt{\eps}}
\left(\frac{x-1}{x+1}\right)^{\Re\beta}
\left(\frac{1-x\eps}{1+x\eps}\right)^{-\Re\beta}
\frac{dx}{x^{1/2}}\cdot
\int_1^{1/\eps}
\left(\frac{y-1}{y+1}\right)^{\Re\beta}
\left|\left(\frac{1-y\eps}{1+y\eps}\right)^{-\beta/2}
-1\right|^2\frac{dy}{y^{3/2}}
$$
and
$$
\int_{1}^{1/\sqrt{\eps}}
\left(\frac{x-1}{x+1}\right)^{\Re\beta}
\left|\left(\frac{1-x\eps}{1+x\eps}\right)^{-\beta/2}
e^{-2Rx\eps}-1\right|^2\frac{dx}{x^{1/2}}\cdot
\int_1^{1/\eps}
\left(\frac{y-1}{y+1}\right)^{\Re\beta}
\frac{dy}{y^{3/2}}.
$$
By a subtitution $x\mapsto 1/x$, $y\mapsto 1/y$, these integrals
become precisely the above integrals (with $\beta$ replaced by
$-\beta$) except that in one integral a term $e^{-2R\eps/x}$
appears, which does affect not the argumentation.
Hence also these terms converge to zero as $\eps\to0$.
\end{proof}

In view of the following proposition, let us make the following
observations. {}From Lemma \ref{l3.19} and Lemma \ref{l3.20}
it follows that 
$$
\det(I\pm K_\beta^0)\iv (I\pm K_{\beta,n})\quad
\mbox{and}\quad
\det(I\pm K_\beta^0)\iv (I\pm \hK_{\beta,R})
$$
are well defined operator determinants for
$-3/2<\Re\beta<1/2$ (in the ``$+$''-case) and
$-1/2<\Re\beta<1$ (in the ``$-$''-case).
Moreover, by Lemma \ref{l3.21} these operator determinants are 
nonzero for sufficiently large $n$ and $R$.

Furthermore it follows that the operator determinant 
$$
\det (I\pm H_\beta^0)\iv(I\pm H_\beta)
$$
is well defined for $-1<\Re\beta<1/2$ (``$+$''-case)
and $-1/2<\Re\beta<3/2$ (``$-$''-case), respectively.
This operator determinant represents a not identically 
vanishing analytic function in $\beta$ (since it equals $1$ for $\beta=0$), and thus it 
is nonzero except possibly on a discrete set.

Finally, from Lemma \ref{l3.19} and its proof we can conclude that 
the determinants $\det(I\pm \Pe K_{\beta}^0\Pe)$ are nonzero for all $\beta$
satisfying $-3/2<\pm\Re\beta<1/2$ and all $\eps>0$.

\begin{proposition}\label{p3.23}
Let $-1<\pm\Re\beta<1/2$. Then for  each $n\ge1$ and $R>0$ we have
\bqn
\lim_{\eps\to0}
\frac{\det(I\pm K_{\beta,\eps,n})}{\det(I\pm \Pe K_\beta^0\Pe)}
&=& \det (I\pm H_\beta^0)\iv(I\pm H_\beta)\cdot
\det(I\pm K_\beta^0)\iv (I\pm K_{\beta,n}),\qquad
\label{l3.74}
\\[1ex]
\lim_{\eps\to0}
\frac{\det(I\pm \hK_{\beta,\eps,R})}{\det(I\pm \Pe K_\beta^0\Pe)}
&=& \det (I\pm H_\beta^0)\iv(I\pm H_\beta)\cdot
\det(I\pm K_\beta^0)\iv (I\pm \hK_{\beta,R}).
\label{l3.75}
\eqn
\end{proposition}
\begin{proof}
First of all we can write
$$
\det(I\pm \hat{K}_{\beta,\eps,R}) =
\det(I\pm \Pe K_\beta^0\Pe)
\det(I\pm A_\eps\pm B_\eps),
$$
where
\bqn
A_\eps &=&  (I\pm \Pe K_\beta^0\Pe)\iv \Per(\hat{K}_{\beta,\eps,R}-
\Pe K_\beta^0\Pe),
\nn\\[1ex]
B_\eps &=&  (I\pm \Pe K_\beta^0\Pe)\iv \Pr (\hat{K}_{\beta,\eps,R}-
\Pe K_\beta^0\Pe).
\nn
\eqn
Equation (\ref{f.as2}) along with the  fact that $\hK_{\beta,R}-K_\beta^0$ 
is trace class implies that
$$
A_\eps = A+o_1(1),\qquad A:=(I\pm K_\beta^0)\iv(\hK_{\beta,R}-K_\beta^0)
$$
(see also Lemma \ref{l3.19}). Similarly, equation (\ref{f.as4}) implies
$$
Y_\eps B_\eps Y_\eps^* = B+o_1(1),\qquad B:=(I\pm H_\beta^0)\iv (H_\beta-H_\beta^0).
$$
Moreover,
$$
B_\eps A_\eps= Y_\eps^*B Y_\eps A+o_1(1)=o_1(1)
$$
since $Y_*\to 0$ weakly.
Hence we can conclude that
\bqn
\lim_{\eps\to 0}
\frac{\det(I\pm \hat{K}_{\beta,\eps,R})}{\det(I\pm \Pe K_\beta^0\Pe)}
&=&
\lim_{\eps\to 0}
\det(I\pm A_\eps)\,\det(I\pm B_\eps)\nn\\
&=& \det(I\pm A)\,\det (I\pm B),\nn
\eqn
which proves the assertion (\ref{l3.75}).
The case of the determinant $\det(I\pm K_{\beta,\eps,n})$
can be treated analogously.
\end{proof}

The previous proposition puts us in position to identify the limit on the
right hand side of (\ref{f3.P}).

\begin{proposition}\label{p3.24}
Let $-3/2<\Re\beta<1/2$ (``$+$''-case) or $-1/2<\Re\beta<1$ (``$-$''-case),
respectively. Then all sufficiently large $n$ and $R$, 
\bqn
\frac{\det P_R(I\pm H(\hat{u}_\beta))\iv P_R}
{\det P_n(I\pm H(u_\beta))\iv P_n}
&=&
\det (I\pm K_{\beta,n})\iv(I\pm \hat{K}_{\beta,R}).
\label{f3.78}
\eqn
\end{proposition}
\begin{proof}
For $-1<\pm \Re\beta<1/2$ and $\beta$ not belonging to a certain
discrete set (namely the set where $\det (I\pm H_\beta^0)\iv (I\pm H_\beta)$ is zero),
we can take the quotient of (\ref{l3.74}) and (\ref{l3.75})
and obtain 
\bqn
\lim_{\eps\to0}
\frac{\det(I\pm \hK_{\beta,\eps,R})}{\det(I\pm K_{\beta,\eps,n})}
&=&
\det (I\pm K_{\beta,n})\iv(I\pm \hat{K}_{\beta,R}).
\eqn
Recall that in a remark made after Proposition \ref{p2.16} we have observed that
$\det(I\pm K_{\beta,\eps,n})$ is nonzero for $\eps>0$ sufficiently small. 
Applying Proposition \ref{p2.16} we obtain identity (\ref{f3.78})
under the assumptions that $-1<\pm\Re\beta<1/2$ and that
$\beta$ does not belong to a certain discrete subset.
We can remove this extra assumption since both sides of the equality are analytic
in $\beta$.
\end{proof}

\begin{theorem}\label{t3.25}
Let $-3/2<\Re\beta<1/2$ (``$+$''-case) or $-1/2<\Re\beta<1$ (``$-$''-case),
respectively.
Then 
\bqn
\det P_R(I\pm H(\hat{u}_\beta))\iv P_R &\sim&
\det P_n(I\pm H(u_\beta))\iv P_n
\eqn
as $R,n\to\iy$ and $R=2n+O(1)$.
\end{theorem}
\begin{proof}
This follows from the previous proposition in connection with Lemma \ref{l3.20}
and Lemma \ref{l3.21}.
\end{proof}


\subsection{Proof of the main results and remarks}
\label{s3.6}

Now are able to prove the main results.

\begin{proofof}{Theorem \ref{t.main}}
We notice first that the proof of the first statement in Theorem
\ref{t.main}(b) follows from Proposition \ref{p2.13}(b).

{}From Theorem \ref{t.Ddiscrete} and Proposition \ref{p2.7}
it follows that 
\bqn
\det\Big[ P_n(I+H(u_\beta))\iv P_n\Big] &\sim&
n^{\beta^2/2+\beta/2} (2\pi)^{-\beta/2} 2^{-\beta^2/2} 
\frac{G(1/2)}{G(1/2-\beta)}, \quad n\to\iy,
\eqn
for $-3/2<\Re\beta<1/2$ and
\bqn
\det\Big[ P_n(I-H(u_\beta))\iv P_n\Big] &\sim&
n^{\beta^2/2-\beta/2} (2\pi)^{-\beta/2} 2^{-\beta^2/2} 
\frac{G(3/2)}{G(3/2-\beta)}, \quad n\to\iy,
\eqn
for $-1/2<\Re\beta<3/2$. With $n=[R/2]$ we can apply Theorem \ref{t3.25},
and we obtain
\bqn
\det\Big[ P_R(I+H(\hu_\beta))\iv P_R\Big] &\sim&
R^{\beta^2/2+\beta/2} (2\pi)^{-\beta/2} 2^{-\beta^2-\beta/2} 
\frac{G(1/2)}{G(1/2-\beta)}, \quad R\to\iy,
\eqn
for $-3/2<\Re\beta<1/2$ and
\bqn
\det\Big[ P_R(I-H(\hu_\beta))\iv P_R \Big]&\sim&
R^{\beta^2/2-\beta/2} (2\pi)^{-\beta/2} 2^{-\beta^2+\beta/2} 
\frac{G(3/2)}{G(3/2-\beta)}, \quad R\to\iy,
\eqn
for $-1/2<\Re\beta<1$.
The proof can now be completed by applying Proposition
\ref{p2.13}.
\end{proofof}

Let us conclude with some final observation. The 
results of the previous sections allow us to 
establish formulas for the determinants of 
$\det P_n(I\pm H(u_{\beta}))\iv P_n$ and
$\det P_R(I\pm H(\hu_{\beta}))\iv P_R$
in terms of certain operator determinants.
We are able to evaluate some (but not all) of these 
determinants explicitly.

The formulas that we obtain might give rise to an alternative
(perhaps clearer) proof of the main result in the sense that one
avoids taking the quotient of the determinants corresponding to the discrete and continuous part
right from the beginning. This would eliminate the annoying discussion of the
non-vanishing of several determinants.

Before establishing the formulas for 
$\det P_n(I\pm H(u_{\beta}))\iv P_n$ and
$\det P_R(I\pm H(\hu_{\beta}))\iv P_R$, we are going to evaluate the asymptotics
of a truncated Wiener-Hopf determinant with a specific, well-behaved symbol.
The result might be of interest in its own since precisely this symbol appear also 
elsewhere. 

\begin{lemma}\label{l3.26}
Let $\phi_\beta(\xi)=1-\sin(\pi\beta)\sech(\pi\xi)$, $\xi\in\R$, and $-3/2<\Re\beta<1/2$. Then
$$
\det W_s(\phi_\beta)\sim e^{-s(\beta/2+\beta^2/2)} 
\frac{G^2(3/2+\beta/2)G^2(1+\beta/2)G^2(1-\beta/2)G^2(1/2-\beta/2)}
{G(1/2)G(3/2)G(3/2+\beta)G(1/2-\beta)}
$$
as $s\to\iy$.
\end{lemma}
\begin{proof}
Using the Achiezer-Kac formula (see e.g. \cite[Sect.~10.80]{BS}), we obtain
$$
\det W_{s}(\phi_\beta)\sim G[\phi_\beta]^{s}E[\phi_\beta], \quad s\to\iy,
$$
where $G[\phi_\beta]$ is given by (\ref{f.55}) and evaluates to
$\exp(- \beta/2-\beta^2/2)$. The constant $E[\phi_\beta]$ is given by
\bqn
E[\phi_\beta] &=&
\exp\left(\int_{0}^\iy x \,(\cF(\log \phi_\beta)(x))
(\cF(\log \phi_\beta)(-x))\,dx
\right)\nn\\
&=&
\exp\left(   \frac{-i}{2\pi}\int_{-\iy}^\iy
(\log\phi_{\beta,+})'(x)(\log\phi_{\beta,-})(x)\,dx
\right).\nn
\eqn
where $\cF$ is the Fourier transform (\ref{f.FT}). The functions 
$\phi_{\beta,\pm}$ stand for the factors of the Wiener-Hopf factorization 
$\phi_\beta$. It is possible to compute these factors explicitly, and one
obtains $\phi_{\beta,\pm}(x)=\psi_\beta(\mp i x/2)$ with
$$
\psi_{\beta}(z) =
\frac{\Gamma(3/4+z)\Gamma(1/4+z)}{\Gamma(3/4+\beta/2+z)
\Gamma(1/4-\beta/2+z)}.
$$
This function is analytic in the right-half plane and has the appropriate behavior at infinity.
Replacing $\phi_{\beta,\pm}$ by $\psi_\beta$ and making a change of variables $z=ix$ gives
$$
E[\phi_\beta]=\exp\left(
\frac{1}{2\pi i}\int_{+i\iy}^{-i\iy}
\frac{\psi_\beta'(-z)}{\psi_\beta(-z)}\, \log\psi_\beta(z)\,dz\right)
$$
A complex function argument implies that this equals the exponential of the residues
of the expression under the integral in the right half plane. Notice that due to the logarithmic 
derivative only simple poles are involved. Thus $E[\phi_\beta]$ equals the exponential of 
the sum ($n=0,1,\dots$) of
$$
\log\psi_\beta(n+3/4)+\log\psi_\beta(n+1/4)-\log\psi_\beta(n+3/4+\beta/2)
-\log\psi_\beta(n+1/4-\beta/2)
$$
A straightforward computation now gives
\bqn
E[\phi_\beta] &=& \prod_{n=0}^\iy
\frac{\Gamma(3/2+n)\Gamma(1+n)}{\Gamma(3/2+\beta/2+n)\Gamma(1-\beta/2+n)}\cdot
\frac{\Gamma(1+n)\Gamma(1/2+n)}{\Gamma(1+\beta/2+n)\Gamma(1/2-\beta/2+n)}
\nn\\[1ex]
&&\times\;
\frac{\Gamma(3/2+\beta+n)\Gamma(1+n)}{\Gamma(3/2+\beta/2+n)\Gamma(1+\beta/2+n)}\cdot
\frac{\Gamma(1+n)\Gamma(1/2-\beta+n)}{\Gamma(1-\beta/2+n)\Gamma(1/2-\beta/2+n)}.
\nn
\eqn
Using the recursion relation for the Barnes $G$-function we obtain that
\bqn
E[\phi_\beta]&=&
\frac{G^2(3/2+\beta/2)G^2(1+\beta/2)G^2(1-\beta/2)G^2(1/2-\beta/2)}
{G(1/2)G(3/2)G(3/2+\beta)G(1/2-\beta)}\cdot R\nn
\eqn
where
\bqn
R &=& \lim_{n\to\iy}\frac{G(1/2+n)G(3/2+n)G(3/2+\beta+n)G(1/2-\beta+n)G^4(1+n)}
{G^2(3/2+\beta/2+n)G^2(1+\beta/2+n)G^2(1-\beta/2+n)G^2(1/2-\beta/2+n)}.\nn
\eqn
Using (\ref{f.BarAsy}) we conclude that $R=1$, which settles the assertion.
\end{proof}

\begin{theorem}
Let $-1<\pm \Re\beta<1/2$. Then for all $n\ge 1$ and $R>0$ we have
\bqn
\det P_n(I\pm H(u_{\beta}))\iv P_n
&=& C_{\pm \beta }\cdot
\det (I\pm H_\beta^0)\iv(I\pm H_\beta)\cdot
\det(I\pm K_\beta^0)\iv (I\pm K_{\beta,n}),\qquad
\\[1ex]
\det P_R(I\pm H(u_{\beta}))\iv P_R
&=&C_{\pm \beta} \cdot\det (I\pm H_\beta^0)\iv(I\pm H_\beta)\cdot
\det(I\pm K_\beta^0)\iv (I\pm \hK_{\beta,R}),
\eqn
where $$C_{\beta}= 2^{\beta^{2}} \frac{G(1/2)G(3/2)G(3/2+\beta)G(1/2-\beta)}
{G^2(3/2+\beta/2)G^2(1+\beta/2)G^2(1-\beta/2)G^2(1/2-\beta/2)}.$$
\end{theorem}
\begin{proof}
We obtain these fromulas from 
the identities  (\ref{f3.Pn}) and (\ref{f3.PR}), from
Propositions \ref{p2.6} and \ref{p2.10} and from
Proposition \ref{p3.23} 
with the constants
$$
C_{\pm\beta} =\lim_{\eps\to0}
\frac{\det(I\pm H(\hu_{\beta,\eps}))}{\det(I\pm \Pi_{[\eps,1]} K_\beta^0 \Pi_{[\eps,1]})}.
$$
Notice in this connection that
$\det(I\pm H(\hu_{\beta,\eps}))=\det(I\pm H(u_{\beta,r}))$ for
$\eps=\frac{1-r}{1+r}$. It remains to evaluate these constants $C_{\pm\beta}$.
This will be done in two steps by establishing an asymptotic formula for
$\det(I\pm H(u_{\beta,r}))$ as $r\to 1$ and for
$\det (I \pm \Pi_{[\eps,1]} K_\beta^0 \Pi_{[\eps,1]})$ as $\eps\to 0$.

For the evaluation of $\det(I\pm H(u_{\beta,r}))$ we rely on the results
of \cite{BE0} (see Theorem 2.5 and formulas (1.12) and (2.15) therein). These results say
that for a sufficiently smooth  nonvanishing function $b$ on $\T$ 
with winding number zero the identity
$$
\det(I+T\iv(b)H(b))=\left(\frac{b_+(1)}{b_+(-1)}\right)^{1/2}\exp\left(-\frac{1}{2}
\sum_{k=1}^\iy k[\log b]_k^2\right)
$$
holds, where $b_+$ is the plus-factor of the Wiener-Hopf factorization of $b$.
We apply this formula with $b(t)=b_+(t)=(1-rt)^{\beta}$ and
$b(t)=b_+(t)=(1+rt)^{\beta}$, respectively. We notice that 
$$\det(I+T\iv(b)H(b))=\det(I+H(b_+)T(b_+\iv))=\det(I+H(b_+\tilde{b}_+\iv)),$$
which is equal to $\det(I\pm H(u_{\beta,r}))$. Notice that we rely on formula (\ref{f.Wflip})
in the ``$-$''-case.
The evaluation of the right hand side gives
$$
\det(I\pm H(u_{\beta,r}))=\left(\frac{1-r}{1+r}\right)^{\pm \beta/2}(1-r^2)^{\beta^2/2}\sim 
\eps^{\pm \beta/2+\beta^2/2}\,2^{\beta^2}, \quad \eps\to0.
$$

The determinant $\det(I\pm \Pe K_\beta^0 \Pe)$ can be expressed as the determinant of a
finite trunctation of a Wiener-Hopf operator.
Proceeding as in the proof of Lemma \ref{l3.19}, the operator 
$K_\beta^0$ is unitarily equivalent to a Wiener-Hopf operator $W(a_\beta)$ with the 
symbol $a_\beta(\xi)=-\sin(\pi\beta)\sech(\pi\xi)$, $\xi\in \R$, while
the projections $\Pe$ transform into $\Pi_{[0,-\log\eps]}$. Thus
$$
\det(I\pm \Pe K_\beta^0 \Pe)=\det W_{-\log\eps}(1\pm a_\beta)
$$
Applying Lemma \ref{l3.26} with $s=-\log\eps$ and $\phi_{\pm\beta}=1\pm a_\beta$
we obtain that $\det(I\pm \Pe K_\beta^0 \Pe)$ is asymptotically equal to
$\eps^{\pm\beta/2+\beta^2/2}$ times the product of the Barnes functions
(with $\beta$ replaced by $-\beta$ in the ``$-$''-case) appearing in Lemma \ref{l3.26}.

Combining both asymptotics yield the desired expression of $C_{\pm\beta}$.
\end{proof}


\end{document}